\documentclass[UTF-8,reqno]{amsart}
\usepackage{enumerate}
\usepackage{mhequ}
\usepackage{amsthm,amsmath,amssymb,url,color, booktabs,nccmath}
\usepackage[left=3.2cm,right=3.2cm,top=4cm,bottom=4cm]{geometry}
\usepackage{mathrsfs}
\usepackage{enumitem,dsfont}
\usepackage{subfigure}
\usepackage[graphicx]{realboxes}
\usepackage{tikz}
\usepackage{color}
\usepackage[colorlinks=true]{hyperref}
\hypersetup{
    linkcolor=blue,          
    citecolor=red,        
    filecolor=blue,      
    urlcolor=cyan
}

\definecolor{darkergreen}{rgb}{0.0, 0.5, 0.0}

\numberwithin{equation}{section}
\def\theequation{\arabic{section}.\arabic{equation}}
\newcommand{\be}{\begin{eqnarray}}
\newcommand{\ee}{\end{eqnarray}}
\newcommand{\ce}{\begin{eqnarray*}}
\newcommand{\de}{\end{eqnarray*}}
\newtheorem{theorem}{Theorem}[section]
\newtheorem{lemma}[theorem]{Lemma}
\newtheorem{proposition}[theorem]{Proposition}
\newtheorem{Examples}[theorem]{Example}
\newtheorem{corollary}[theorem]{Corollary}

\newtheorem{definition}[theorem]{Definition}
\theoremstyle{definition}
\newtheorem{remark}[theorem]{Remark}

\DeclareMathOperator{\supp}{supp}

\def\[{{\Big[}}
\def\]{{\Big]}}
\def\<{{\langle}}
\def\>{{\rangle}}
\def\({{\Big(}}
\def\){{\Big)}}

\def\bx{{\mathbf{x}}}

\def\dif{{\mathord{{\rm d}}}}

\def\min{{\mathord{{\rm min}}}}

\def\={&\!\!=\!\!&}

\newcommand{\norm}[1]{\left\|#1\right\|}

\def\cF{{\mathcal F}}

\def\cP{{\mathcal P}}

\def\mN{{\mathbb N}}

\def\mP{{\mathbb P}}

\def\mR{{\mathbb R}}

\def\mT{{\mathbb T}}

\def\mZ{{\mathbb Z}}

\def\bP{{\mathbf P}}
\def\bQ{{\mathbf Q}}

\def\1{{\mathbf{1}}}

\def\geq{\geqslant}
\def\leq{\leqslant}

\def\le{\leqslant}

\def\div{\mathord{{\rm div}}}

\def\[{{\Big[}}
\def\]{{\Big]}}
\def\<{{\langle}}
\def\>{{\rangle}}
\def\({{\Big(}}
\def\){{\Big)}}

\def\bx{{\mathbf{x}}}

\def\dif{{\mathord{{\rm d}}}}

\def\min{{\mathord{{\rm min}}}}

\def\={&\!\!=\!\!&}
\def\bt{\begin{theorem}}
\def\et{\end{theorem}}
\def\bl{\begin{lemma}}
\def\el{\end{lemma}}
\def\br{\begin{remark}}
\def\er{\end{remark}}
\def\bx{\begin{Examples}}
\def\ex{\end{Examples}}
\def\bd{\begin{definition}}
\def\ed{\end{definition}}
\def\bp{\begin{proposition}}
\def\ep{\end{proposition}}
\def\bc{\begin{corollary}}
\def\ec{\end{corollary}}

\def\geq{\geqslant}
\def\leq{\leqslant}

\def\le{\leqslant}

\def\div{\mathord{{\rm div}}}

\def\bP{{\mathbf P}}

 \def\R{\mathbb R}
 \def\R{\mathbb R}    
\def\N{\mathbb N}  
   
\def\<{\langle} \def\>{\rangle}

\allowdisplaybreaks

\begin{document}

\title[Sharp Non-uniqueness of SDEs on $\mR^d$]{Sharp Non-uniqueness  in Law  for Stochastic Differential Equations on the Whole Space}

\author{Huaxiang L\"u}
\address[H. L\"u]{Academy of Mathematics and Systems Science, Chinese Academy of Sciences, Beijing 100190, China}
\email{lvhuaxiang22@mails.ucas.ac.cn }

\author{Michael R\"ockner}
\address[M. R\"ockner]{ Fakult\"at f\"ur Mathematik, Bielefeld Universit\"at, D 33615 Bielefeld, Germany, and  Academy of Mathematics and Systems Science, Chinese Academy of Sciences, Beijing 100190, China}
\email{roeckner@math.uni-bielefeld.de}
 

\thanks{
Funded by the Deutsche Forschungsgemeinschaft (DFG,German Research Foundation)-Project-ID 317210226-SFB 1283"}

\begin{abstract}
In this paper, we investigate the stochastic differential equation  on $\mathbb{R}^d,d\geq2$:
\begin{align*}
    \dif X_t&=v(t,X_t)\dif t+\sqrt{2} \dif W_t.
\end{align*}
For any finite collection of initial probability measures $\{\mu^i_0\}_{1\leq i\leq M}$ on $\mathbb{R}^d$ and $\frac{d}{p}+\frac{1}{r}>1$, we construct a divergence-free drift field $v\in L_t^rL^p\cap C_tL^{d-}$ such that the associated SDE admits at least two distinct weak solutions originating from each initial measure $\mu^i_0$. This result is sharp in view of the well-known uniqueness of strong solutions for drifts in $C_tL^{d+}$, as established in \cite{KR05}.
As a corollary,  there exists a measurable set $A\subset\mathbb{R}^d$ with positive Lebesgue measure such that for any $x\in A$, the SDE with drift $v$ admits at least two weak solutions when  with start in $x\in A$. 
The proof proceeds by constructing two distinct probability solutions to the associated Fokker-Planck equation via a convex integration method  adapted to   all of  $\mathbb{R}^d$ (instead of merely the torus),  together with  refined heat kernel estimate.  

\end{abstract}


\date{\today}

\maketitle

\tableofcontents

\section{Introduction}
In this paper, we are concerned with  stochastic differential equations (SDEs)  on the whole space $\R^d,d\geq2$, of type
\begin{align}
    \dif X_t&=v(t,X_t)\dif t+\sqrt{2}\dif W_t,\ \ t\in[0,T], 
    \label{eq:sde}\\
 \mathcal{L}_{ X_0}&=\mu_0,\notag
\end{align}
where $T>0$, $v : [0 , T] \times  \mR^d \to \mR^d$ is a Borel function,   $W_t$ is a standard $\mR^d$-valued Brownian motion on some probability space $(\Omega, \mathcal{F}, \mathbf{P})$, and $\mu_0$ is a probability measure on $\mR^d$.

The well-posedness of SDEs has been a central topic in both pure and applied mathematics, due to their fundamental role in modeling dynamical systems with random perturbations. To analyze the criticality of equation \eqref{eq:sde}, we employ the following scaling transformation: $$W^\epsilon_t=\epsilon^{-1}W_{\epsilon^2t},\ X^\epsilon_t=\epsilon^{-1}X_{\epsilon^2t},\ v^\epsilon(t,x)=\epsilon v(\epsilon^2t,\epsilon x),$$   which preserve the law of the Brownian motion. Under this scaling, the behavior of the rescaled drift $v^\epsilon$ in the Lebesgue space $L^r_tL^p:=L^r([0,T];L^p(\mathbb{R}^d))$ is characterized by:
\begin{align*}
    \|v^\epsilon\|_{L^r_tL^p}=\epsilon^{1-\frac{d}{p}-\frac{2}{r}}\|v\|_{L^r_tL^p}.
\end{align*}
As $\epsilon\to0$,  SDE \eqref{eq:sde} is classified as subcritical when 
$\frac dp+\frac2r<1$, in which case the associated quantity vanishes. The critical regime corresponds to $\frac dp+\frac2r=1$ and 
  the supercritical regime corresponds to $\frac dp+\frac2r>1$, where the quantity diverges, indicating that the Brownian noise is no longer sufficient 
to counteract the singularity of the drift.

For $p=r=\infty$, Veretennikov \cite{Ver80} was the first to demonstrate the uniqueness of strong solutions using Zvonkin’s transform \cite{Zvo74}.   In the subcritical case,   Krylov and the second named author \cite{KR05} established the existence and uniqueness of strong solutions to \eqref{eq:sde}. Then the well-posedness problem with general multiplicative noise was studied in the foundational work \cite{Zha05,Zha11}. Subsequent developments are contained in \cite{BC03,FGP10,Zha11,FF13,Zha16,BFGM19,XXZZ20,RZ21a}. 

 The critical regime  presents  subtle challenges.  Beck, Flandoli, Gubinelli and Maurelli \cite{BFGM19} proved that for $\|v\|_{L^d}$ small enough,  the SDE admits a strong solution starting from a diffusive random variable in a certain class.  A significant progress was made by Krylov \cite{Kry20} who proved the strong
well-posedness  of SDEs in the case  $ v\in L^d (\mR^d )$. Then, for the case $d \geq 3$, $v\in  L^r_tL^p,d<p\leq \infty$, or $v \in C_t L^d$, or $v\in L^\infty_tL^d$ with divergence-free condition, Zhao and the second named author \cite{RZ23} proved that  \eqref{eq:sde} admits a unique weak solution  within a class satisfying some Krylov-type estimate.   We refer to \cite{Nam20,Kry20b,Kry20c,RZ21,RZZ25} for further results.

However, in the supercritical case, the well-posedness problem is
much more challenging and remains not fully understood. On the one hand, under the additional assumption that the drift is divergence-free,  weak well-posedness has been established in   certain classes: Zhang and Zhao \cite{ZZ21} established weak existence and uniqueness in the sense of approximation for  $\frac{d}{p} + \frac{2}{r} \leq 2$;    Galeati \cite{Gal24}  demonstrated strong existence and uniqueness  for  \eqref{eq:sde}, under the assumption that $v$ is a Leray solution to the 3D Navier-Stokes equations obtained through approximation.    We also refer to \cite{Gal23,HZ23,GP24,Gra24} for more results in the supercritical regime. On the other hand, notable counterexamples demonstrating non-uniqueness in law for \eqref{eq:sde} have also been constructed:
  For dimensions $d \geq 1$ and $p > d$, Galeati and Gerencs\'er \cite{GG25} constructed an example (allowed to be not divergence-free), where non-uniqueness in law holds when starting from the origin.  In $d \geq 3$, the second named author, Zhang and Zhao \cite{RZZ25} exhibited a divergence-free drift in weak Lorentz space $L^{p,\infty},$ $\frac d2<p<d$, that admits two distinct weak solutions starting from the origin. When considering the torus case,   sharp non-uniqueness in law in dimension 2 was shown in \cite{LRZ25}.

In particular, if we restrict to the case $r=\infty$, then for any vector field  $v\in L_t^\infty L^p$, $p>d$, \eqref{eq:sde} admits a unique strong solution (see \cite{KR05}). For the case $d\geq3$, for any $v\in C_t L^d$, or divergence-free $v\in L_t^\infty L^d$, the SDE admits a unique weak solution within a class such that some Krylov-type estimate holds (see \cite{RZ23}). However, for a given initial probability measure on $\mR^d$, the question of whether uniqueness in law holds for SDEs with drifts $v\in L_t^\infty L^p$ in the supercritical regime $p<d$ 
 remains open, especially for the two-dimensional  case. Notably, in all known examples in the literature,  non-uniqueness in law   occurs only when  starting from the origin. 
 
It is also worth mentioning that, in the supercritical case, the divergence-free condition on the drift  matters. Actually, if the drift  is not divergence-free, there  is an example in \cite[Example 7.4]{BFGM19} such that the SDE may not have weak solutions   in $L^p(\mR^d),p<d$. On the other hand, the  divergence-free condition 
helps when considering the existence and conditional weak uniqueness,  cf. \cite{ZZ21,HZ23,RZZ25}.
In summary, it is particularly intriguing to investigate divergence-free drift terms that may lead to a failure of weak uniqueness. 

\subsection{Main result}
Regarding the aforementioned discussion, in this paper, we make a new   contribution to the understanding of weak (non-)uniqueness in  SDEs by constructing  a divergence-free drift $v\in C_tL^p,p<d$ that leads to non-uniqueness in law, which covers the optimal range for the weak well-/ill-posedness of the SDE.

Let $\mathcal{P}$ denote the space of probability measures on $\mathbb{R}^d$. For any finite collection of initial probability measures $\{\mu^i_0\}_{1\leq i\leq M} \subset \mathcal{P}$ and  $\frac{d}{p} + \frac{1}{r} > 1$, our main result is to construct a divergence-free drift $v \in L^r_t L^p$  such that  non-uniqueness in law holds for the respective   SDE with drift $v$ started from each of these initial measures: \begin{align}
    \dif X^i_t&=v(t,X^i_t)\dif t+\sqrt{2}\dif W_t,
    \label{eq:sdei}\\
   \mathcal{L}_{X^i_0}&=\mu^i_0.\notag
\end{align} 

\bt\label{thm:4converge}
  Let $d\geq2,\gamma\in(0,1),1<s<d$ and $p,r\in[1,\infty]$ satisfying $\frac dp+\frac1r>1$. For a finite collection of initial distributions $\{\mu^i_0\}_{1\leq i\leq M}\subset\cP$, there exists   a divergence-free
vector field $v\in L^r([0,T];L^p) \cap C([0,T];L^s)
$ such that each SDE  \eqref{eq:sdei} admits at least two distinct weak solutions $X^{i,1},X^{i,2}$.
Moreover, for $1\leq i\leq M,j=1,2$,   it holds that 
\begin{itemize}
    \item $\mathbf{E}[\int_0^T|v(s,X^{i,j}_s)|^{1+\epsilon}\dif s] <\infty$ for some $\epsilon>0$;
    \item for $t>0$, the solutions admit densities $\rho^{i,j}$  satisfying $\rho^{i,j}(t)\geq (1-\gamma)e^{t\Delta}\mu^i_0$. 
\end{itemize} 

\et

Our work also extends previous results  in \cite{GG25,RZZ25}, which    were limited to initial distributions concentrated at the origin. We demonstrate that weak non-uniqueness can occur for arbitrary initial probability measures. Through  decomposition of probability measures, we further establish that this non-uniqueness in law holds for a set of initial positions $x$ with positive Lebesgue measure, thereby substantially broadening the scope of possible non-uniqueness scenarios.

\begin{corollary}\label{thm:nonuni:pointwise}
Let $d\geq2,1<s<d$ and $p,r\in[1,\infty]$ satisfying $\frac dp+\frac1r>1$.  Then there exists   a divergence-free
vector field $v\in L^r([0,T];L^p) \cap C([0,T];L^s)
$ and a measurable  set $A\subset \mathbb{R}^d$ with positive Lebesgue measure such that  for every $x \in A$, there are at least two distinct weak solutions to \eqref{eq:sde} starting at $\delta_x$. 
\end{corollary}

Unlike the previous counterexamples, in this paper we establish a stronger form of non-uniqueness, specifically, the non-uniqueness in the evolution of particle densities. These densities solve the associated Fokker-Planck equations, which read as
\begin{align}\label{eq:fpe:i}
\partial_t \rho^i-\Delta\rho^i+\div(v\rho^i)&=0,\ \ \div v=0,\\
\rho^i(0)\dif x&=\mu^i_0(\dif x),\notag
\end{align}
where $\{\mu^i_0\}_{1\leq i\leq M}$ is a collection of probability measures on $\mR^d$. 
A solution to \eqref{eq:fpe:i} is meant in the following weak sense:
\bd \label{def:sol:fpe}
A non-negative function $\rho\in C((0,T];L^1)$ is called a weak solution to \eqref{eq:fpe:i} if 
$$\int_0^T\int_{\mR^d}|v(t,x)|\rho(t,x)\dif x\dif t<\infty,$$
and for every test function  $f\in C_c^2(\mR^d)$, and every $t\in(0,T]$
\begin{align}
    \int_{\mR^d}f(x)\rho(t,x)\dif x-\int f(x)\mu_0^i(\dif x)=\int_0^t\int_{\mR^d}\rho(s,x)(\Delta f+v\cdot\nabla f)(s,x)\dif x\dif s\label{def:solution}
.\end{align}
\ed

With the help of the superposition principle introduced in \cite[Theorem 2.5]{Tre16},  Theorem \ref{thm:4converge} is a direct result of the following non-uniqueness result for the Fokker-Planck equation \eqref{eq:fpe:i}, which we shall establish  through the convex integration method,  together with  refined heat kernel estimates.

\bt\label{thm:non_pde_4}
 Let $d\geq2,\gamma\in(0,1),1<s<d$ and $p,r\in[1,\infty]$ satisfying $\frac dp+\frac1r>1$.  Then there exists   a divergence-free
vector field $v\in L^r([0,T];L^p) \cap C([0,T];L^s)
$, such that each Fokker-Planck equation \eqref{eq:fpe:i} admits at least two distinct solutions $\rho^{i,1},\rho^{i,2}\in C((0,T];L^1)
$ in the sense of Definition \ref{eq:fpe:i}.  

Moreover, for $1\leq i\leq M,j=1,2$,  it holds that 
 $\int_{[0,T]}\int|v|^{1+\epsilon}\rho^{i,j}\dif x\dif s <\infty$ for  some $\epsilon>0$
    and that   $\rho^{i,j}(t)\geq (1-\gamma)e^{t\Delta}\mu^i_0$ for $t>0$.

\et

Finally, we  note that  by linearity, \eqref{eq:fpe:i}  admits infinitely many distinct solutions for every $1\leq i\leq M$. Consequently, the associated SDE \eqref{eq:sdei} also has  infinitely many distinct weak solutions.

The main result will be established  using the convex integration technique. 
This technique was first introduced to fluid dynamics by De Lellis and Sz\'ekelyhidi Jr. \cite{DLS09, DLS10, DLS13}. This method has  led to numerous groundbreaking results for fluid dynamics on the torus. For the incompressible Euler equations, the famous Onsager conjecture was proved in \cite{Ise18,BDLSV19}. For the Navier-Stokes equations,  sharp non-uniqueness of weak solutions has been shown in \cite{BVb,BCV18,CL22,CL23}.   We refer to \cite{TCP,DSJ17, GKN23,NV23,GKN24,GR24,BCK24b}  for more progress on the Euler or Navier-Stokes  equations. Recently, this method has been applied to fluid dynamics on the whole space, see \cite{MNY24a,MNY24b}.
For more details and references, interested readers are referred to the comprehensive reviews \cite{BV,BV21}. 

Concerning   transport equations  on the tours,  we refer to \cite{CGSW15,MS18,MS19,MS20,BCDL21,CL21,PS23,CL22b,MS24} for recent developments concerning non-uniqueness in transport equations and related ODEs.
We also note that the convex integration method has been successfully applied to   stochastic fluid dynamics, see \cite{BFH20,HZZ20,Yam22, HZZ,HZZ21b,HLP24,HZZ19,Pap24,MS24,HZZ22,LZ23b,LZ24,LRZ25}.

\subsection{Idea of the proof}
 As mentioned earlier, we will prove the main result using the convex integration method. This method is primarily developed in the torus setting, and for the convenience of readers who may not be so familiar with it, in Appendix~\ref{app:a}  we provide an overview of how the convex integration method works on the torus.

Our work presents a novel extension of the convex integration method to establish non-uniqueness of positive solutions for Fokker-Planck equations on the whole space $\mR^d$, going beyond previous periodic domain results.  The transition from the torus to unbounded domains introduces substantial new challenges  that require new probabilistic and analytical arguments. These technical innovations entail, to our knowledge, the first successful application of convex integration methods to prove non-uniqueness for positive solutions in non-periodic settings. 

\begin{enumerate}
   \item 
When working in the whole space setting, maintaining compact support for the stress error remains crucial for implementing convex integration. However, a fundamental difficulty arises in the iteration process: we cannot use  the inverse of the divergence operator directly, since even for a function with compact support, the inverse of its divergence  may still have support in the whole space. To overcome this obstacle, we introduce a new decomposition of the stress term (detailed in Section \ref{tamr}) that separates it into two distinct components - a principal part, expressed as the divergence of a compactly supported function, and a residual global error term with small bounded $C_{t,x}^0$ norm. This decomposition enables us to proceed with the convex integration scheme by handling the principal part through standard techniques, while simultaneously controlling the global error via a specially constructed global perturbation.  

 \item  When constructing global perturbations, we face the fundamental challenge of preserving   positivity of the solution, which differs significantly from those encountered in fluid equations. In addition, we do not have an invariant measure as a reference measure in contract to the torus case. The core of our approach lies in using the heat-propagated measure $e^{t\Delta}\mu_0$   as a reference solution, around which we build two distinct positive solutions through carefully designed perturbations. A critical technical requirement emerges from the spatial decay of $e^{t\Delta}\mu_0$: our perturbations must precisely match its exponential decay at infinity to maintain global positivity. This is achieved through a novel strategy by constructing drift terms with compact support, coupled with a careful analysis of the influence in exterior regions. The key innovation centers on our treatment of the global perturbation $\theta^{\rm (glo)}$ generated by the drift, which solves the heat equation:
\begin{align*}
    \partial_t\theta^{\rm (glo)}-\Delta\theta^{\rm (glo)}=F,
\end{align*} 
  with a specially designed forcing term 
$F$. Here $F$ is constructed with compact spatial support and vanishes near the time origin $t=0$, enabling the crucial proper exponential decay \begin{align}\label{key:est}
  \bigg| \int_{0}^te^{(t-s)\Delta}F(s)\dif s \bigg|\ll e^{t\Delta}\mu_0 
\end{align}
  through Proposition \ref{prop:heatkernel}   using   heat kernel estimates. The required smallness condition is attained via an additional iteration of the inverse divergence operator developed in Section \ref{tamr}, combined with the periodic properties of our building blocks. 
  
 We further note that the key estimate \eqref{key:est} appears to be new and, to the best of our knowledge, has not been addressed in previous literature. This estimate plays a crucial role in our analysis, and its possible generalization to more general diffusion operators remains an open and particularly interesting problem.

 \item 
To achieve the supercritical regime $\frac dp+\frac1r>1$, we introduce the $L^m$-based intermittent spatial-time jets, incorporating an additional degree of freedom through the parameter $m\in[1,\infty]$. This additional  freedom stems from a  feature of transport-type equations, contrasting  with the $L^2$-based constraint in fluid equations.  

\end{enumerate}
Beyond those already mentioned, the proof also requires several additional modifications for the $\mathbb{R}^d$ setting, such as  a carefully tuned mollification procedure  and a refined construction of the principal density perturbations. We will detail these technical adaptations in the   proof in the main text of this proof.

\noindent{\bf Organization of the paper.}  First,  in Section \ref{gij}, we introduce the $L^m$-based  intermittent spatial-time jets, which play a crucial role in the convex integration construction. Subsequently, in Section \ref{cogpss3}, we establish Theorem~\ref{thm:4converge} through a proof of the non-uniqueness of solutions to the Fokker-Planck equations, as stated in Theorem \ref{thm:non_pde_4}. The proof of the main convex integration iteration procedure is presented in Section \ref{proof:propcase4}.   In Appendix \ref{app:a}, we give an overview of the convex integration method on the torus.

\noindent{\bf Notations.} Let $T>0$, $\mN_{0}:=\mN\cup \{0\}$. Throughout the manuscript, we write  $\mT^d = \mR^d/\mZ^d$ for the $d$-dimensional flat torus, and identify $\mathbb{T}^d$-valued functions with  periodic functions on $\mathbb{R}^d$.
\begin{itemize}
    \item 
We employ the notation $a\lesssim b$ if there exists a constant $c>0$ such that $a\leq cb$.  
\item   Given a Banach space $E$ with a norm $\|\cdot\|_E$, we write $C_tE=C([0,T];E)$ for the space of continuous functions from $[0,T]$ to $E$, equipped with the supremum norm.  For $p\in [1,\infty]$ we write $L^p_tE=L^p([0,T];E)$ for the space of $L^p$-integrable functions from $[0,T]$ to $E$, equipped with the usual $L^p$-norm. We also wirte $L^p_{[a,b]}E=L^p([a,b];E)$.

\item  
For $\alpha\in(0,1)$ we  define $C^\alpha_tE$ to be the space of $\alpha$-H\"{o}lder continuous functions from $[0,T]$ to $E$, endowed with the norm $\|f\|_{C^\alpha_tE}=\sup_{s,t\in[0,T],s\neq t}\frac{\|f(s)-f(t)\|_E}{|t-s|^\alpha}+\sup_{t\in[0,T]}\|f(t)\|_{E},$   and write  $C_t^\alpha$ in the case when $E=\mathbb{R}$. 
 
\item  We use $L^p$ to denote the set of  standard $L^p$-integrable functions.
For $s>0$, $p>1$ we set $W^{s,p}:=\{f\in L^p; \|(I-\Delta)^{\frac{s}{2}}f\|_{L^p}<\infty\}$ with the norm  $\|f\|_{W^{s,p}}=\|(I-\Delta)^{\frac{s}{2}}f\|_{L^p}$.

\item For $N\in \N_0 $, $C^N$ denotes the space of $N$-times differentiable functions equipped with the norm
	$$
	\|f\|_{C^N}:=\sum_{\substack{|\alpha|\leq N, \alpha\in\N^{d}_{0} }}\| D^\alpha f\|_{ L^\infty_x}.$$ Similarly, if the norm is taken in space-time, we use the notation $C^{N}_{t,x}$.

\item For any $\mT^d$-valued smooth function $f$, we define  the projections $\mathbb{P}_{=0} f := \int_{\mathbb{T}^d} f\dif x$, and  $\mathbb{P}_{\neq0} f := f -\int_{\mathbb{T}^d} f\dif x$. 

\item By $p(t,x)=(4\pi t)^{-\frac d2}e^{-|x|^2/4t}$ we denote the $d$-dimensional heat kernel for $\Delta$ on $\mR^d$.

\item  We denote the Lebesgue measure on $\mR^d$ by $\dif x$.
\item  We denote the law of a random  variable $X$ by $\mathcal{L}_X$.
\end{itemize}

\section{$L^m$-based  intermittent spatial-time jets}\label{gij}

In this section we introduce the the notion of $L^m$-based intermittent spatial-time jets, where  $m\in[1,\infty]$, which can be seen as a generalization of the  $L^2$-based building blocks  presented in  \cite[Appendix C.1]{LRZ25}. In this section, the building blocks are defined on the torus and are regarded as periodic functions on $\mathbb{R}^d$.

First, we introduce  the following geometrical lemma:
\bl\cite[Lemma 3.1]{BCDL21}\label{lem:cv2}
Let $d\geq2$. There exists a finite set $\Lambda\in \mathbb{S}^{d-1}\cap \mathbb{Q}^d $ and non-negative $C^\infty$-function  $\Gamma_{\xi}:\mathbb{S}^{d-1}\to\R $ such
that for every $R \in \mathbb{S}^{d-1}$
$$R=\sum_{\xi\in\Lambda}\Gamma_\xi(R)\xi.$$
\el
 With Lemma \ref{lem:cv2} in hand,   it is easy to generate $2$ disjoint families $\Lambda^{1},\Lambda^{2}$, where each one enjoys the property of Lemma \ref{lem:cv2} by taking suitable rational rotations of one fixed set. For simplicity, we denote $\Lambda:= \Lambda^{1}\cup\Lambda^{2}$. Moreover, we know that  $\{\Gamma_{\xi}\}_{\xi\in \Lambda }$ are uniformly bounded.

For parameters $\lambda,r_\perp, r_\parallel > 0 $, we assume
$$\lambda^{-1}\ll r_{\perp}\ll r_{\parallel}\ll 1,\ \ \lambda r_{\perp}\in\mathbb{N}.$$

For each $\xi\in\Lambda$ let us define $A^i_\xi\in \mathbb{S}^{d-1}\cap \mathbb{Q}^d,\ i=1,2,...,d-1$, such that $\{\xi, A^i_\xi,i=1,...,d-1\}$
form an orthonormal basis in $\mathbb{R}^d$. 
 Let  $n_*\in\mN$ such that$\{n_*\xi, n_*A^i_\xi,i=1,...,d-1\}\subset\mathbb{Z}^d$
for every $\xi\in\Lambda$. 

We define $\phi : \mathbb{R}^{d-1} \to \mathbb{R}$ to be a smooth function supported in the unit ball, 
such that $\phi \equiv 1$ on $B(0,\tfrac{1}{3})$ and has zero mean. 
We then define $\Phi$ by $\phi = -\Delta \Phi$. 
The existence of such functions can be found in \cite[Appendix~C.1]{LRZ25}. 
Let $\psi : \mathbb{R}\to\mathbb{R}$  be a smooth, mean-zero function with
support in $B(0,1)$ satisfying $\psi\equiv 1$ on $B(0,\frac13)$. 

Define $\phi' : \mathbb{R}^{d-1} \to \mathbb{R}$ to be a smooth non-negative function with support in $B(0,\frac13)$ satisfying $$\int_{\mathbb{R}^{d-1}}\phi'(x_1,x_2,...,x_{d-1})\dif x_1\dif x_2..\dif x_{d-1}=1,$$
and let $\psi' : \mathbb{R} \to \mathbb{R}$ be a smooth non-negative function with support in $B(0,\frac13)$ such that
$$\int_\mathbb{R}\psi'(x_d)\dif x_d=1.$$
Then, it is straightforward to verify that
\begin{align}
    \phi\phi'=\phi',\ \ \psi\psi' =\psi'.\label{phi2}
\end{align}

Let  $m\geq1$ be fixed, we define the rescaled cut-off functions
$$\phi _{r_{\perp},m}(x_1,x_2,...,x_{d-1})=\frac1{r_{\perp}^{(d-1)/m}}\phi (\frac{x_1}{r_{\perp}},\frac{x_2}{r_{\perp}},...,\frac{x_{d-1}}{r_{\perp}}),$$
$$\Phi _{r_{\perp},m}(x_1,x_2,...,x_{d-1})=\frac1{r_{\perp}^{(d-1)/m}}\Phi (\frac{x_1}{r_{\perp}},\frac{x_2}{r_{\perp}},...,\frac{x_{d-1}}{r_{\perp}}),$$
$$\psi _{r_{\parallel},m}(x_d)=\frac1{r_{\parallel}^{1/m}}\psi (\frac{x_d}{r_{\parallel}}).$$
Similarly,  for a conjugate exponent $m'\in[1,\infty]$ satisfying $\frac1{m'}+\frac1m=1$, we define $\phi '_{r_{\perp},m'}$ and $\psi'_{r_{\parallel},m'}$ as the same manner. Then we periodize $\phi _{r_{\perp},m},\Phi _{r_{\perp},m},$ $\psi _{r_{\parallel},m}$, $\phi '_{r_{\perp},m'}$ and $\psi'_{r_{\parallel},m'}$ so that they can be viewed as functions on $\mathbb{T}^{d-1}$ and $\mathbb{T}$ respectively.
Consider a large time oscillation parameter $\mu > 0$. For every $\xi\in\Lambda$ we introduce
$$\psi _{(\xi,m)}(t,x):=
\psi _{r_{\parallel},m}(n_*r_{\perp}\lambda(x\cdot \xi-\mu t)),$$
$$\Phi _{(\xi,m)}(x):=
\Phi _{r_{\perp},m}(n_*r_{\perp}\lambda x\cdot A^1_{\xi},...,n_*r_{\perp}\lambda x\cdot A^{d-1}_{\xi}),$$
$$\phi _{(\xi,m)}(x):=
\phi _{r_{\perp},m}(n_*r_{\perp}\lambda x\cdot A^1_{\xi},...,n_*r_{\perp}\lambda x\cdot A^{d-1}_{\xi}).$$
Similarly  we define  $\phi '_{(\xi,m')}$  and $\psi'_{(\xi,m')}$. 

The building blokes $W_{(\xi,m)} :\mathbb{R} \times\mathbb{T}^d \to \mathbb{R}^d$ and $\Theta_{(\xi,m')} :\mathbb{R} \times\mathbb{T}^d \to \mathbb{R}$ are defined as
$$W_{(\xi,m)}(t,x):=\xi\psi_{(\xi,m)}(t,x)\phi_{(\xi,m)}(x),$$
$$\Theta_{(\xi,m')}(t,x):=\psi'_{(\xi,m')}(t,x)\phi'_{(\xi,m')}(x).$$
By the definition and \eqref{phi2} we have that
\begin{align}
\int _{\mathbb{T}^d}W_{(\xi,m)}\Theta_{(\xi,m')}\dif x=\xi,\label{eq:intwthe}\\
\partial_t\Theta_{(\xi,m')}+\mu r_\perp^{\frac{d-1}{m}}r_\parallel^{\frac1m}\div (W_{(\xi,m)}\Theta_{(\xi,m')})=0.\label{eq:ptthe+}
\end{align}

Since $W_{(\xi,m)}$ is not divergence-free, we
introduce the skew-symmetric  corrector term
\begin{align}
V_{(\xi,m)}:=\frac{1}{(n_*\lambda)^2}(\xi\otimes\nabla\Phi_{(\xi,m)}-\nabla\Phi_{(\xi,m)}\otimes\xi)\psi_{(\xi,m)}\notag
\end{align}
satisfying
\begin{align}
\div V_{(\xi,m)}
=W_{(\xi,m)} -\frac{1}{(n_*\lambda)^2}\nabla\Phi_{(\xi,m)}\xi\cdot\nabla\psi_{(\xi,m)}.\label{divOmega}
\end{align}
Finally, we obtain that for $N, M \geq 0$ and $p\in [1, \infty]$ the
following holds
\begin{align}
\|\nabla^N\partial_t^M\psi_{(\xi,m)}\|_{C_tL^p(\mT^d)}\lesssim r_\parallel^{\frac{1}{p}-\frac{1}{m}}(\frac{r_\perp\lambda}{r_\parallel})^N(\frac{r_\perp\lambda\mu}{r_\parallel})^M,
&\label{int2}\\
\|\nabla^N\phi_{(\xi,m)}\|_{L^p(\mT^d)}+\|\nabla^N\Phi_{(\xi,m)}\|_{L^p(\mT^d)}\lesssim r_\perp^{\frac{d-1}{p}-\frac{d-1}m}\lambda^{N},
&\label{int3}\\
\|\nabla^N\partial_t^MW_{(\xi,m)}\|_{C_tL^p(\mT^d)}+\lambda\|\nabla^N\partial_t^MV_{(\xi,m)}\|_{C_tL^p(\mT^d)}
\lesssim r_\perp^{\frac{d-1}{p}-\frac{d-1}m}r_\parallel^{\frac{1}{p}-\frac{1}{m}}\lambda^{N}(\frac{r_\perp\lambda\mu}{r_\parallel})^M,&\label{int4}\\ \|\nabla^N\partial_t^M\Theta_{(\xi,m')}\|_{C_tL^p(\mT^d)}\lesssim r_\perp^{\frac{d-1}{p}-\frac{d-1}{m'}}r_\parallel^{\frac{1}{p}-\frac{1}{m'}}\lambda^{N}(\frac{r_\perp\lambda\mu}{r_\parallel})^M,& \label{int4theta}
\end{align}
where the implicit constants may depend on $p,m,N$ and $M$, but are independent of $\lambda,r_\perp,r_\parallel,\mu$. These estimates can be easily deduced from the definitions.

Then let us  introduce a family of temporal functions   to oscillate the building blocks intermittently in time.
Let  $K\in\mN$ be fixed, and $G \in C_c^\infty(0, 1)$ be non-negative and 
$$\ \ \int_0^1G^2(t)\dif t=1.$$
  Let $\eta >0$  be a small constant satisfying $ \eta K\ll1$. For $\xi\in\Lambda$ as defined above,  and $1\le i\le K,$ we define $\tilde{g}_{(\xi,i,m)}: \mathbb{T}\to \mathbb{R}$ as the 1-periodic extension of $\eta^{-1/m}G(\frac {t-t_{(\xi,i)}}{\eta})$,  where $t_{(\xi,i)}$ are chosen so that $\tilde{g}_{(\xi,i,m)}$ have disjoint supports for distinct $(\xi,i)$.  
We will also oscillate the perturbations at a large
frequency $\sigma\in\mathbb{N}$. So, we define
$$g_{(\xi,i,m)}(t)=\tilde{g}_{(\xi,i,m)}(\sigma t).$$

For the corrector term we define $H_{(\xi,i,m)},h_{(\xi,i,m)}:\mathbb{T}\to\mathbb{R}$ by
\begin{align}
H_{(\xi,i,m)}(t)=\int_0^{t}g_{(\xi,i,m)}(s)\dif s,\ \ h_{(\xi,i,m)}(t)=\int_0^{\sigma t}(\tilde{g}_{(\xi,i,m)}(s)\tilde{g}_{(\xi,i,m')}(s)-1)\dif s, \label{eq:parth}
\end{align}
where we recall that $\frac1{m'}+\frac1m=1$, In view of the zero-mean condition for $\tilde{g}_{(\xi,i,m)}(t)\tilde{g}_{(\xi,i,m')}(t)-1$, we see that $h_{(\xi,i,m)}$ is $\mathbb{T}/\sigma$-periodic, and for any $N\geq0,p\geq1$
\begin{align}
    \|g_{(\xi,i,m)}\|_{W^{N,p}}\lesssim(\frac{\sigma}{\eta})^N\eta^{1/p-1/m},\  \  \|h_{(\xi,i,m)}\|_{L^\infty}\leq1,\label{bd:gwnp}
\end{align}
where the universal constant is independent of the choices of  $i$ and $\xi$.

\section{Proof of Main result}\label{cogpss3}
In this section, we prove our main result, Theorem~\ref{thm:4converge}. 
Without loss of generality, we set $T = 1$ from now on. 
The proof begins by establishing Theorem~\ref{thm:non_pde_4} at the PDE level, 
where we construct non-unique solutions to the Fokker-Planck equation~\eqref{eq:fpe:i}. 
To this end, we consider a system of $N_0 = 2M$ equations given by
\begin{align}\label{5:eq:fpe+in}
\partial_t \rho^i-\Delta\rho^i+\div(v \rho^{in,i})+\div(v \rho^i)&=0,\ \
\div v=0,\\
\rho^i(0)=0,\notag
\end{align}
where the initial densities are defined for $1 \leq k \leq M$ by
$$\rho^{in,2k-1}(t)=\rho^{in,2k}(t):=e^{t\Delta} \mu_0^i.$$  
In what follows, the drift $v$ will be constructed so that it vanishes near the origin $t = 0$. 
It is then straightforward to verify that $\rho^{ in,i} + \rho^i$ satisfies 
equation~\eqref{eq:fpe:i}.

Regarding equation~\eqref{5:eq:fpe+in}, for parameters satisfying 
$1 < s < d$ and $\frac{d}{p} + \tfrac{1}{r} > 1$, 
our goal is to prove Theorem~\ref{thm:non_pde_4} by constructing a divergence-free drift
 $ v \in L^r_tL^p\cap L_{t}^{d_0}L^{d_0}\cap C_tL^s$
such that for each $i$, equation~\eqref{5:eq:fpe+in} admits a solution
$\rho^i\in L_{t}^{d_0'}L^{d_0'}\cap C_tL^1$,
where $d_0$ and $d_0'$ are conjugate exponents satisfying 
$\frac{1}{d_0} + \frac{1}{d_0'} = 1$. 
Moreover, we show that $\rho^{2k-1} \neq \rho^{2k}$ for $1 \leq k \leq M$, 
thereby demonstrating the non-uniqueness of solutions.

 The proof employs a convex integration scheme specifically adapted to construct positive solutions 
to the Fokker--Planck equations on the whole space~$\mathbb{R}^d$. 
In Section~\ref{sec:estrhoin}, we establish several key estimates for the initial components 
$\rho^{in,i}$. 
Subsequently, the framework of the convex integration iteration and the corresponding iterative procedure 
are developed in Proposition~\ref{prop:case4} within Section~\ref{sec:setup}. 
The proof of Theorem~\ref{thm:non_pde_4} is then presented in Section~\ref{sec:proofmaireuslt}, 
followed by the derivation of Theorem~\ref{thm:4converge} through an application of the superposition principle.

\subsection{Estimate on $\rho^{in,i}$}\label{sec:estrhoin}
We begin by deriving several fundamental estimates for the initial densities~$\rho^{{in},i}$,
which will play a crucial role in the subsequent analysis.
Recall that $p(t,x)$ denotes the $d$-dimensional heat kernel associated with the Laplacian~$\Delta$. We define the averaged heat kernel
\begin{align*}
    \overline p(t,x):= \int_{[-\frac12,\frac12]^d}p(t,x-y)\dif y.
\end{align*} 

By definition and the positivity of the heat kernel $p(t,x)$, it is straightforward to verify that $\rho^{in,i}(\frac1{12},x)>0$ for any $1\leq i\leq N_0,x\in\mR^d$. Moreover, since $ [-1/2,1/2]^d$ is a compact domain, there exists $c_{in}>0$ such that for any $1\leq i\leq N_0,x\in [-1/2,1/2]^d$,
\begin{align*}
\rho^{in,i}(\frac1{12},x)\geq c_{in}.
\end{align*}
 Then, it holds that for $t\in[\frac1{12},1]$,
\begin{align}
\rho^{in,i}(t,x)&= \int_{\mR^d}p(t-\frac1{12},x-y)\rho^{in,i}(\frac1{12},y)\dif y\notag\\ 
&\geq c_{in}\int_{[-\frac12,\frac12]^d}p(t-\frac1{12},x-y)\dif y=c_{in}\overline p(t-\frac1{12},x).\notag
\end{align}
In particular,   for $(t,x)\in[\frac16,1]\times  [-1/2,1/2]^d$, by definition, it holds that 
\begin{align*}
    \overline p(t-\frac1{12},x)\geq \int_{[-\frac12,\frac12]^d}(4\pi)^{-d/2}e^{-3d}\dif y=:c_d,
\end{align*}
which implies that for $(t,x)\in[\frac16,1]\times  [-1/2,1/2]^d$, it holds   that
\begin{align}
    \rho^{in,i}(t,x)\geq c_{in}\overline p(t-\frac1{12},x)\geq c_{in}c_d.\label{bd:rhoingeqcd}
\end{align}
 
Moreover, by the smoothing property of the heat kernel, 
there exists a constant $0 < C_{in} < \infty$ such that, 
for all $1 \le i \le N_0$, the following estimate holds:

\begin{align}
\|\rho^{in,i}\|_{C_{[\frac{1}{6},1],x}^1} \leq C_{in}. \label{5:bd:rhoin}
\end{align}

\subsection{Convex integration set up}\label{sec:setup}  

We now apply the convex integration iteration to the system~\eqref{5:eq:fpe+in}, 
with the iteration indexed by $q \in \mathbb{N}_0$.
We consider an increasing sequence $\{\lambda_q\}_{q\in\mathbb{N}_0}\subset\mathbb{N}$ which diverges to $\infty$, and a sequence $\{\delta_q\}_{q\in\mathbb{N}_0}\subset (0,1]$
 which is decreasing to 0. We choose $a\in\mathbb{N}_0,\beta\in(0,1)$ and $b\in  \mathbb{N}$. Let
$$\lambda_q=a^{(b^q)},q\geq0,\ \ \delta_q=\epsilon_0^{d+1}\lambda_1^{2\beta}\lambda_q^{-2\beta},q\geq1,\ \ \delta_0=1,$$
where $\beta>0$ will be chosen sufficiently small and $a,b$ will be chosen sufficiently large. Here $\epsilon_0\in(0,\frac1{12}]$ is a small universal constant to be determined later. Moreover, we use the estimate $$\sum_{q\geq1} \delta_q^{1/(d+1)}\leq \epsilon_0\sum_{q\geq1}a^{(1-q)2b\beta/(d+1)}\leq \frac{\epsilon_0}{1-a^{-2b\beta/(d+1)}}< 2\epsilon_0,$$
which boils down to
\begin{align}\label{ieq:4ab2}
a^{2b\beta/(d+1)}>2,
\end{align}
assumed from now on.

At each step $q$, a pair $(v_q,\rho^{i}_q,M^{i}_q) _{1\leq i\leq N_0}$ is constructed solving the following system on $[0,1]$:
\begin{align}\label{eq:4qth}
\partial_t \rho^{i}_q-\Delta\rho^{i}_{q}+\div(v_q \rho^{in,i})+\div(v_q \rho^{i}_q)&=-\div M^{i}_q,\ \
\div v_q=0,
\end{align}
 where 
 $M_q^i$ are some vector fields.

To handle the initial condition, we require that $\rho_q^i=0$  on $[0,T_q]$, where 
$$   T_q := \frac13 - \sum_{ 1\leq r\leq q} {\delta}_r^{1/2}. $$
By applying  \eqref{ieq:4ab2} we obtain $\frac16< T_q\leq\frac13$. Here and in the following we define $\sum_{1\leq r\leq 0}:=0.$  To control the spatial support of the drift and the stress term, we further define
\begin{align*}
     \Omega_q := \[-\frac13-\sum_{ 1\leq r\leq q} {\delta}_r^{1/2},\frac13+\sum_{ 1\leq r\leq q} {\delta}_r^{1/2}\]^d\subset [-\frac12,\frac12]^d.
\end{align*}

With the above assumptions in place, the main iteration scheme is formulated as follows:

\bp\label{prop:case4}
Under the assumption of Theorem \ref{thm:4converge}, there exist $d+1>d_0>2>d_0'>1$ with $\frac{1}{d_0}+\frac{1}{d_0'}=1$ and a choice of parameters $a,b,\beta$  such
that the following holds: Let $(v_q,\rho^{i}_q,  M^{i}_q) _{1\leq i\leq N_0}$ be a solution to the system \eqref{eq:4qth} satisfying  $\int_{\mR^d}\rho^{i}_q\dif x=0$,
\begin{align}\label{bd:4vql2}
\|v_q\|_{L^{d_0}_tL^{d_0}}\leq C_vC_0^{1/d_0}\sum_{m=0}^q\delta_{m}^{1/d_0},\ \ \|\rho^{i}_q\|_{L^{d_0'}_tL^{d_0'}}\leq C_\rho  C_0^{1/{d_0'}}\sum_{m=0}^{q}\delta_m^{1/{d_0'}}
\end{align}
for some universal constants $C_0,C_v, C_\rho\geq1$, and
\begin{align}
\|v_q\|_{C_{t,x}^1}\leq C_0^{1/d_0}\lambda_q^{d+4},\
\|\rho^{i}_{q}\|_{C_{t,x}^1}
\leq C_0^{1/d_0'} \lambda_q^{d+4},
\label{bd:4rhoqc1}\\
 \|M^{i}_q\|_{L^1_tL^1}\leq C_0 \delta_{q+1},\ \ \|\partial_t M^{i}_{q}\|_{L_t^1L^1}+  \| \nabla M^{i}_{q}\|_{L_t^1L^1}\leq \lambda_{q}^{2d+8},&\label{bd:4rql1}\\
 \rho^{i}_q(t)=v_q(t)=M^{i}_q(t)=0\ {\rm on}\ [0,T_q],&\label{bd:4mq=rhoq=0}\\
 \supp v_q,\supp M^{i}_q\subset \Omega_q.\label{bd:supp}
\end{align}
Then there exists 
$(v_{q+1},\rho^{i}_{q+1},  M^{i}_{q+1})_{1\leq i\leq N_0} $ which solves \eqref{eq:4qth} and satisfies \eqref{bd:4vql2}-\eqref{bd:supp} at the level $q+1$ and
\begin{align}
\|v_{q+1}-v_{q}\|_{L^{d_0}_tL^{d_0}}\leq
C_vC_0^{1/d_0}\delta_{q+1}^{1/d_0},\label{bd:4vq+1-vqldd}\ \
\|\rho^{i}_{q+1}-\rho^{i}_{q}\|_{L^{d_0'}_tL^{d_0'}}\leq
C_\rho C_0^{1/{d_0'}}\delta_{q+1}^{1/{d_0'}}.
\end{align}
Moreover, we have 
\begin{align}
\|v_{q+1}-v_{q}\|_{L^{r}_tL^{p}}\leq
\delta_{q+1}^{1/d_0},\ \ \|v_{q+1}-v_q\|_{C_tL^s}\leq
\delta_{q+1}^{1/d_0},\label{bd:4vq+1-vqlpr}\\
\|\rho^{i}_{q+1}-\rho^{i}_{q}\|_{C_tL^1}\leq
  \delta_{q+1}^{1/d_0'},\label{bd:4rhoq+1-rhoql1}\\
(\rho^{i}_{q+1} - \rho^{i}_q)(t,x)\geq- \delta_{q+1}^{1/d_0'} \overline p(t-\frac1{12},x)\ {\rm for\ }t\in[\frac16,1].\label{bd:4rhoq+1-rhoqgeq}
\end{align}
\ep
Here $C_0$ is determined by the choice of the starting iterations, and $C_v,C_\rho$ are two constants determined by the improved H\"older's inequality for $v_q$, $\{\rho_q^i\}_{1\leq i\leq N_0} $ respectively, and other implicit constants in the proof. 
We further note that the term $v_q\rho^{in,i}$ is always well-defined by noticing \eqref{5:bd:rhoin}, \eqref{bd:4rhoqc1} and the fact that $v_q=0$ on $[0,T_q]$. 

Proposition~\ref{prop:case4} constitutes the central technical component of this paper. 
A detailed proof of the proposition will be presented in Section~\ref{proof:propcase4} below.

\subsection{Proof of main result}\label{sec:proofmaireuslt}
In this section, we first present the proof of Theorem~\ref{thm:non_pde_4}, 
assuming the validity of Proposition~\ref{prop:case4}.

\begin{proof}[Proof of Theorem \ref{thm:non_pde_4}]
Let $\gamma\in(0,1)$ be fixed.
Let $F(t,x)$ be a smooth (not divergence-free) bounded $\mR^d$-valued function with support in $[\frac13,1]\times \Omega_0$, and satisfying $$\|\div F\|_{C_{t,x}^0}\leq \frac\gamma 2 c_{in}c_d.$$ 
Then  we denote  $$ \|\div F\|_{C_tL^1}=: c_{F}>0.$$

Let $\epsilon_0=\min\{\frac{\gamma c_{in}}{4},\frac{c_F}2,\frac1{12}\}$. 
 We intend to start the iteration from 
$(v_0, \rho^{i}_0,  M^{i}_0)_{1\leq i\leq N_0}$ which are defined as 
$$\rho^{i}_0=(-1)^i \div F,\ \  v_0=0,\ \ M^{i}_0=(-1)^{i+1} (\partial_t-\Delta) F.$$
By choosing $C_0$ large enough (depending on $\epsilon_0$), we have 
\begin{align}
   \| \rho^{i}_0\|_{L^{d_0'}_tL^{d_0'}}+ \|\rho^{i}_0\|_{C_{t,x}^1}  \leq  C_0^{1/d_0'},\ \ \|M^{i}_0\|_{L^1_tL^1}\leq \epsilon_0^{d+1} C_0.\notag
\end{align}
 Then \eqref{bd:4vql2}-\eqref{bd:supp} are satisfied as $\delta_0=1,\delta_1=\epsilon_0^{d+1}$.

Moreover, by \eqref{bd:rhoingeqcd} we have  
\begin{align}
    \frac \gamma2\rho^{in,i}+\rho^i_0\geq(\frac \gamma2\rho^{in,i}+\rho^i_0 )1_{[\frac16,1]\times [-1/2,1/2]^d}\geq \frac\gamma2c_{in}c_d-\frac\gamma2c_{in}c_d=0.\label{eq:rhoin+rhoi0}
\end{align}

Next, we use Proposition \ref{prop:case4} to build inductively $(v_q,\rho^{i}_q, M^{i}_q) $ for every $q \geq 1$. By \eqref{ieq:4ab2} and \eqref{bd:4vq+1-vqldd}-\eqref{bd:4rhoq+1-rhoql1}, the sequence $\{(v_q,\rho^{i}_q)\}_{q\in \N}$ is
Cauchy in 
\begin{align*}
    \(L^r([0,1];L^p)\cap L^{d_0}([0,1]\times \mR^d)\cap C([0,1];L^s)\)\times \(L^{d_0'}([0,1]\times \mR^d)\cap C([0,1];L^1)\)^{N_0}\end{align*}
     and we denote by $(v,\rho^{i})$ its limit, where $v$ is also divergence-free with compact support. Since $\int \rho^{i}_q\dif x=0$, we deduce that $\int \rho^{i}\dif x=0$. 
 Clearly by \eqref{bd:4rql1} and \eqref{bd:4mq=rhoq=0},  $\rho^{i}$  solves   \eqref{5:eq:fpe+in}.
 
 Now, we define $\overline \rho^i:=\rho^i + \rho^{{in},i}$, which satisfies \eqref{eq:fpe:i} in the sense of Definition \ref{def:sol:fpe}. Then we  verify that each $\overline \rho^i$ is nonnegative  and $\overline \rho^{2k-1}$ and $\overline \rho^{2k}$ are distinct from one another with the same initial distribution.
In fact, by \eqref{bd:rhoingeqcd}, \eqref{ieq:4ab2}, \eqref{bd:4rhoq+1-rhoql1}-\eqref{eq:rhoin+rhoi0} we have
\begin{align*}
\rho^i+\gamma\rho^{in,i}& \geq  \frac\gamma2\rho^{in,i} +
\sum_{q=0}^\infty (\rho_{q+1}^i- \rho_q^i) \geq (\frac\gamma2c_{in}-\sum_{q=0}^\infty \delta_{q+1}^{1/d_0'})\overline p(t-\frac1{12},x) >0,\ {\rm for}\ t>\frac16,\\
\|\overline \rho^{2k-1}-\overline \rho^{2k}\|_{C_{[\frac23,1]}L^1}&\geq \|\rho^{2k-1}_0-\rho^{2k}_0\|_{C_{[\frac23,1]}L^1}-\|\rho^{2k-1}-\rho^{2k-1}_0\|_{C_{[\frac23,1]}L^1}-\|\rho^{2k}_0-\rho^{2k}\|_{C_{[\frac23,1]}L^1}\notag\\ 
&\geq2( \|\div F\|_{C_tL^1}- \sum_{q=0}^\infty \delta_{q+1}^{1/d_0'})>0.
\end{align*}
Then by the fact that $\rho^i(t)=0$ for $t\in[0,\frac16]$, we imply that $\overline{\rho}^i\geq (1-\gamma)\rho^{in,i}$.

By \eqref{bd:4rhoqc1}, \eqref{bd:4vq+1-vqldd} and interpolation, we have $|v| \in L_{t}^{d_0(1+\epsilon)}L^{d_0(1+\epsilon)}$ for some $\epsilon > 0$ sufficiently small. Since $|\rho^i| \in L_{t}^{d_0'}L^{d_0'}$, we deduce that $|v|^{1+\epsilon} \rho^i \in L_{t}^{1}L^{1}$.  Furthermore, from \eqref{5:bd:rhoin} and the fact that $v$ has compact support satisfying $v = 0$ on $[0,\frac{1}{6}]$, we conclude $|v|^{1+\epsilon} \rho^{in,i} \in  L_t^{1}L^{1}$. Combining these results, we obtain
$|v|^{1+\epsilon} \overline {\rho}^i \in  L_t^{1}L^{1}.$
\end{proof}

Having established Theorem~\ref{thm:non_pde_4}, Theorem~\ref{thm:4converge} follows directly through an application of the superposition principle.
 \begin{proof}[Proof of Theorem \ref{thm:4converge}]
For  $1<s<d,\frac dp+\frac1r>1$, by applying Theorem \ref{thm:non_pde_4}, there are   $v\in L^r([0,1];L^p)\cap C([0, 1]; L^s)$  and a collection of  probability densities $\overline \rho^i(t)$ 
satisfying \eqref{eq:fpe:i}. Then we define $\mu^i_t(\dif x):=\overline \rho^i(t)\dif x$, which forms a family of probability measures. 
Moreover, these measures satisfy
\begin{align*}
     \int_0^1\int|v(s,x)|\mu_s^i(\dif x)\dif s= \int_0^1\int_{\mR^d}|v(s,x)|\overline \rho^i(s,x)\dif x\dif s<\infty,
\end{align*}
and $t \to \mu^i_t$ is weakly continuous on $[0,1]$,
 we are in position to apply the superposition principle (see \cite[Theorem 2.5]{Tre16}) for \begin{align}
\partial_t \mu^i-\Delta\mu^i+\div(v\mu^i)&=0, \ \
\mu^i_{t=0}=\mu^i_0.\notag
\end{align} 
 More precisely, let  $C([0,1];\mR^d) $ be the space of continuous functions,  
 equipped  with its Borel $\sigma$-algebra and its natural filtration generated by the canonical process $\Pi_t, t \in [0,1]$, defined by $$\Pi_t(\omega) := \omega(t),\ \  \omega\in C([0,1];\mR^d).$$
  There exists  a group of  probability measure  $\mathbf{Q}^i$ on $C([0,1];\mR^d) $ which is a martingale solution associated to diffusion operator $$L:=\Delta+v\cdot \nabla.$$ 
  
Then  by a standard result (see \cite[Theorem
 2.6]{Str87}), there exists a $d$-dimensional $(\cF_t)_{t\geq0}$-Brownian motion $W_t$, $t \in [0,1]$,  on a stochastic basis $(\Omega,\cF,(\cF_t)_{t\geq0},\bP)$ and  a group of continuous $(\cF_t)_{t\geq0}$-progressively
 measurable maps $\{X^i_t\}_{1\leq i\leq 2M}: [0,1] \times \Omega \to \mR^d$ satisfying  SDE \eqref{eq:sdei}. Moreover, we have $\bP\circ (X^i_t)^{-1}=\bQ^i\circ \Pi_t^{-1}=\mu^i_t$ for $t\in[0,1]$. Then it is easy to see that $\mathcal{L}_{X^{2k-1}_0}=\mathcal{L}_{X^{2k}_0}=\mu_0^i$, while $\mathcal{L}_{X^{2k-1}_t}\neq \mathcal{L}_{X^{2k}_t}$.
 
\end{proof}

 Corollary \ref{thm:nonuni:pointwise} follows by an argument analogous to the proof of \cite[Theorem 1.7]{LRZ25}.

\section{Proof of Proposition \ref{prop:case4}}\label{proof:propcase4}
In this section, we extend the convex integration method to construct solutions for a system defined on the whole space. The proof follows the  convex integration approach while incorporating necessary adaptations for the unbounded domain. We begin in Section \ref{sec:4choicepara} by determining the choice of parameters. Section \ref{sec:mpll4} then details the mollification procedure. The core construction appears in Section \ref{sec:4defq+1}, where we define the new iteration pair $(v_{q+1},\rho^{i}_{q+1})$  with a carefully designed global perturbation $\theta_{q+1}^{(g,i)}$ that ensures that we can define a suitable  new stress error supported in the domain $\Omega_{q+1}$. Subsequently, Section \ref{sec:4defmq+1} introduces the crucial stress terms $M^{i}_{q+1}$ and the global stress error $F_{q+1}^i$. 
 At this stage, it is essential to apply the inverse divergence operator iteration, as introduced in Section \ref{tamr}, to ensure that the support of $M^{i}_{q+1}$ is located in $\Omega_{q+1}$. This condition is crucial for maintaining the consistency of the iterative process. 
In the end, in Sections \ref{sec:4estwq+1}-\ref{sec:4estmq+1}, we conclude the proof with the verification of all required inductive estimates.  The most technically involved step consists in estimating the global perturbation $\theta_{q+1}^{(g,i)}$ 
  to ensure that it exhibits the decay behavior required by  \eqref{bd:4rhoq+1-rhoqgeq}. This analysis relies crucially on the heat kernel estimates established in Proposition~\ref{prop:heatkernel}.

\subsection{Choice of parameters}\label{sec:4choicepara}
In our analysis, several carefully chosen parameters will play crucial roles, with their values being precisely calibrated to satisfy the intricate network of compatibility conditions required for our estimates.  Let $d\geq2,1<s<d,\frac{d}{p}+\frac{1}{r}>1$ be fixed. First we introduce a integer $N>4d$ large enough  satisfying $$\frac{d}{p}+\frac{1}{r}>1+\frac{4d}N,\ \ \frac{d}{s}>1+\frac{4d}N,$$ and define $$d_0:=d+1-\frac{4d}{N}\in(d,d+1), \ \ d_0':=\frac{d_0}{d_0-1}\in(1,2).$$

Then for the sufficiently small $\alpha\in(0,1)$ to be chosen, we take $ {l}:=\lambda_{q+1}^{-\frac{3\alpha}{2}}\lambda_{q}^{-d-4}$ and have 
 \begin{align}
     {l}^{-1}\leq\lambda_{q+1}^{2\alpha},\ \   l\lambda_q^{2d+8} \leq 
     \lambda_{q+1}^{-\alpha}\ll\delta_{q+2},\ \ \lambda_{q}^{d+4}\ll\lambda_{q+1}^\alpha,\label{para42}
 \end{align} provided $\alpha b> 2d+8,\alpha>2\beta b$. 
   In the sequel, we also need
$(12d+42)\alpha<\frac1{2N}.$
 
 The above can be obtained by choosing $\alpha>0$ small such that $(12d+42)\alpha<\frac1{2N}$, and choosing $b\in 2N\mathbb{N}$ and large enough such that $\alpha b>2d +8$, and finally choosing $0<\beta<\frac{\alpha}{2b}$.

In the end, we  increase $a$ such that \eqref{ieq:4ab2} holds. In the sequel, we also increase $a$  to absorb  various implicit and universal constants in subsequent estimates.

\subsection{Mollification}\label{sec:mpll4}
To avoid the loss of derivative, we first need to mollify the stress term.
Let $\phi_l:=\frac{1}{l^d}\phi(\frac{\cdot}{l})$ be a family of standard radial mollifiers on $\mathbb{R}^d$, and $\varphi_l:=\frac{1}{l}\varphi(\frac{\cdot}{l})$ be a family of standard  mollifiers with support in $(0,1)$.  We define 
\begin{align}
M^{i}_l:
&=(M^{i}_{q}*_x\phi_l)*_t\varphi_l.\notag
\end{align}
For the mollification around $t = 0$,  since  $M_q^i$ vanishes  around $t = 0$, we  can directly extend its value to $t\leq 0$ by $0$. Here we deliberately avoid mollifying the density terms to avoid the  challenge in controlling the difference 
$\rho^i_q-\rho^i_l$, which must maintain exponential decay  crucial to the argument.

Since  $l\leq \delta_{q+1}^{1/2}$, we know that   $M^{i}_l(t)=0$ for $t\in [0,T_{q+1}]$, and $\supp  M^i_l\subset \Omega_{q+1}.$

To  end this section, by the mollification estimates, the spatial-time embedding $W^{d+\frac43,1}\subset L^\infty$  and  \eqref{bd:4rql1} we obtain   for $N\geq 0$,
    \begin{align}
    \|M_l^i\|_{C_{t,x}^N}&\lesssim l^{-d-\frac43-N} \|M_{q }^i\|_{L_t^1L^1}\lesssim C_0  l^{-d-\frac43-N}.\label{bd:4mlcn}
 \end{align}

\subsection{Construction of $v_{q+1}$ and $\rho^{i}_{q+1}$}\label{sec:4defq+1}
In this section, we proceed with the construction on the perturbations on $v_q$ and $\rho_{q}^i$ employing the $L^{d_0}$-based   building blocks introduced in Section \ref{gij}.   First, we define the parameters 
\begin{align}
  \lambda=\lambda_{q+1},\ r_\perp=\lambda^{-1+\frac1{N}},\ r_\parallel=\lambda^{-1+\frac2{N}},\ \eta=\lambda^{-1},\  \mu=r_\perp^{-\frac{d-1}{d_0}}r_\parallel^{-\frac1{d_0}}, \ \sigma=\lambda^{\frac1{2N}}, \label{def:4para}
\end{align}
where we recall that  $d_0,d_0'$  and $N>4d$ are defined in Section \ref{sec:4choicepara}. Then we have
\begin{align}
   r_\perp^{\frac{d-1}{p}-\frac{d-1}{d_0}} r_\parallel^{\frac{1}{p}-\frac{1}{d_0}}\eta^{\frac1r-\frac1{d_0}},\ r_\perp^{\frac{d-1}{s}-\frac{d-1}{d_0}} r_\parallel^{\frac{1}{s}-\frac{1}{d_0}}\eta^{-\frac1{d_0}}\leq\lambda^{-\frac1N},\notag\\
   r_\perp^{d-1-\frac{d-1}{d_0'}} r_\parallel^{1-\frac{1}{d_0'}}\eta^{-\frac{1}{d_0'}},\ \lambda r_\perp^{d-1-\frac{d-1}{d_0'}}r_\parallel^{1-\frac{1}{d_0'}}\eta^{1-\frac1{d_0'}}\leq\lambda^{-\frac1N}.\label{bd:paralam4}
\end{align}
In fact, by a direct calculation we have
    \begin{align*}
        \lambda r_\perp^{d-1-\frac{d-1}{d_0'}}r_\parallel^{1-\frac{1}{d_0'}}\eta^{1-\frac1{d_0'}}=r_\perp^{d-1-\frac{d-1}{d_0'}} r_\parallel^{1-\frac{1}{d_0'}}\eta^{-\frac{1}{d_0'}}\leq r_\parallel^{\frac{d}{d_0}}\eta^{\frac{1}{d_0}-1}= \lambda^{\frac{2d}{Nd_0}+\frac{d_0-d-1}{d_0}}= \lambda^{-\frac{2d}{Nd_0}}\leq \lambda^{-\frac1N},\\
 r_\perp^{\frac{d-1}{p}-\frac{d-1}{d_0}} r_\parallel^{\frac{1}{p}-\frac{1}{d_0}}\eta^{\frac1r-\frac1{d_0}}\leq  \lambda^{(\frac{d+1}{N}-d)(\frac{1}{p}-\frac{1}{d_0})-\frac1r+\frac1{d_0}}
 \leq\lambda^{-\frac{d}{p}-\frac1r+\frac{d+1}{d_0}+\frac{d+1}{N}}\leq \lambda^{-1-\frac{4d}N+1+\frac{2d}{N}+\frac{d+1}{N}}\leq \lambda^{-\frac1N},\\
  r_\perp^{\frac{d-1}{s}-\frac{d-1}{d_0}} r_\parallel^{\frac{1}{s}-\frac{1}{d_0}}\eta^{-\frac1{d_0}}\leq  \lambda^{(\frac{d+1}{N}-d)(\frac{1}{s}-\frac{1}{d_0})+\frac1{d_0}}
 \leq\lambda^{-\frac{d}{s}+\frac{d+1}{d_0}+\frac{d+1}{N}}\leq \lambda^{-1-\frac{4d}N+1+\frac{2d}{N}+\frac{d+1}{N}}\leq \lambda^{-\frac1N}.
    \end{align*}
Here it is required that $b$ is a multiple of $2N$ to ensure that $\lambda_{q+1}r_{\perp}=a^{(b^{q+1})/N}\in\mathbb{N}$ and $\sigma=a^{(b^{q+1})/2N}\in\mN$.
 
Let $\chi\in C_c^\infty(-\frac34,\frac34)$  be a non-negative function such that
$\sum_{n\in\mathbb{Z}} \chi(t - n) = 1$ for every $t \in \R$. Let $\tilde{\chi} \in C_c^\infty(-\frac45,\frac45)$ be a non-negative function satisfying $\tilde{\chi} = 1$ on
$[-\frac34,\frac34]$ and $\sum_{n\in\mathbb{Z}} \tilde{\chi}(t - n) \leq2$.  

Then, we fix a family parameters $ \zeta = 20/\delta_{q+2}$ and consider two disjoint sets $\Lambda^{1}, \Lambda^{2}$ defined in Section \ref{gij}. We use the notation   $ \Lambda^i =  \Lambda^1$ for $i$ odd, and $ \Lambda^i =  \Lambda^2$ for $i$ even. Taking $K = N_0$ in the construction of Section~\ref{gij},  we obtain a family of pairwise disjoint functions $g_{(\xi,i,d_0)}$ and $ g_{(\xi,i,d_0')}$  for $\xi\in\Lambda$ and $1\leq i\leq N_0$.
In the following we apply the building blocks $W_{(\xi,d_0)},\Theta_{(\xi,d_0')}$ and $H_{(\xi,i,d_0)}$ as defined in Section \ref{gij} to define \begin{align*}
    W_{(\xi,n,i)}(x,t):=W_{(\xi,d_0)}\(x,\(\frac{n}{ \zeta}\)^{1/d_0}H_{(\xi,i,d_0)}(t)\),\\
    \Theta_{(\xi,n,i)}(x,t):=\Theta_{(\xi,d_0')}\(x,\(\frac{n}{ \zeta}\)^{1/d_0}H_{(\xi,i,d_0)}(t)\).
\end{align*} 
 Similarly, we define $V_{(\xi,n,i)},\Phi_{(\xi,n,i)}
 $ and all other terms appearing in  Section \ref{gij}. Now by  \eqref{eq:ptthe+}, \eqref{eq:parth} and the choice of $\mu$ we have
 \begin{align}
     \partial_t\Theta_{(\xi,n,i)}+\(\frac{n}{ \zeta}\)^{1/{d_0}}g_{(\xi,i,d_0)}\div (W_{(\xi,n,i)}\Theta_{(\xi,n,i)})=0.\label{eq:4ptthe+n}
 \end{align}

As the next step, we first define the perturbations for the drift term. For every $1\leq i\leq N_0$, we define the principle perturbation and the correct perturbation  as
\begin{align*}
w_{q+1}^{(p,i)}:&=\sum_{n\geq3}\sum_{\xi\in\Lambda^{n}}\tilde{\chi}( \zeta|M^{i}_l|-n)\(\frac{n}{ \zeta}\)^{1/d_0}W_{(\xi,n,i)}g_{(\xi,i,d_0)},\\
w_{q+1}^{(c,i)}:&=\sum_{n\geq3}\sum_{\xi\in\Lambda^{n}}\(-\tilde{\chi}( \zeta|M^{i}_l|-n)\(\frac{n}{ \zeta}\)^{1/{d_0}}\frac{1}{(n_*\lambda_{q+1})^2}\nabla\Phi_{(\xi,n,i)}\xi\cdot\nabla\psi_{(\xi,n,i)}\notag\\ 
&\quad +\nabla( \tilde{\chi}( \zeta|M^{i}_l|-n))\(\frac{n}{ \zeta}\)^{1/{d_0}}:V_{(\xi,n,i)}\)g_{(\xi,i,d_0)}.
\end{align*}
Here we use the notation $(\nabla( \tilde{\chi}( \zeta|M_l^i|-n)):{V}_{({\xi,n,i})})^i:=\sum_{j=1}^d\partial_j( \tilde{\chi}( \zeta|M_l^i|-n)){V}_{({\xi,n,i})}^{ij} $ for $i=1,2,...,d$. 
We remark that here and in the following the first sum runs for $n$ in the range 
\begin{align*}
  3\leq n\leq 1+ \zeta|M^{i}_l|\leq 1+l^{-2/3}\|M^{i}_l\|_{C_{t,x}^0}\leq 1+Cl^{-d-2},
\end{align*}
where we used \eqref{bd:4mlcn} and $C$  is a universal constant.

By \eqref{divOmega} we have $\div ( w_{q+1}^{(p,i)}+w_{q+1}^{(c,i)})=0$ since
\begin{align}
w_{q+1}^{(p,i)}+w_{q+1}^{(c,i)} 
&=\sum_{n\geq3}\sum_{\xi\in\Lambda^{n}}\div\(\tilde{\chi}( \zeta|M^{i}_l|-n)\(\frac{n}{ \zeta}\)^{1/d_0} V_{(\xi,n,i)}\)g_{(\xi,i,d_0)}
,\label{eq:4wq+1p+wq+1c}
\end{align}
and  $V_{(\xi,n,i)}$ is a skew-symmetric matrix.  

Then we define the perturbations for the densities. For $1\leq i\leq N_0$, we define the principle perturbation and the correct perturbation  as
\begin{align}
\theta_{q+1}^{(p,i)}:&=
\sum_{n\geq3}\chi( \zeta|M^{i}_l|-n)\(\frac{n}{ \zeta}\)^{1/{d_0'}}\sum_{\xi\in\Lambda^{n}}\Gamma_{\xi}\(\frac{M^{i}_l}{|M^{i}_l|}\)\(\Theta_{(\xi,n,i)}-\mathbb{P}_{0}\Theta_{(\xi,n,i)}\)g_{(\xi,i,d_0')},\notag\\
\theta_{q+1}^{(c,i)}:&=- \hat{\rho}^i_q\cdot\int_{\mathbb{R}^d}\theta_{q+1}^{(p,i)}\dif x,\label{def:thetaq+1c}
\end{align}
where we recall $\mathbb{P}_{0}f=\int_{\mathbb{T}^d} f\dif x.$ Here   $ \hat{\rho}^i_q$ is some probability density with support in $\Omega_{q+1}$ satisfying $\|\hat{\rho}^i_q\|_{C_{x}^2}\lesssim1.$   

Here, we use of the projected term 
$\mathbb{P}_{\neq 0}\Theta_{(\xi,n,i)}$ rather than $\Theta_{(\xi,n,i)}$ in constructing the principal perturbation. This projection guarantees the applicability of the inverse divergence operation during subsequent stress term construction.

Then by the identity \eqref{eq:intwthe}, Lemma \ref{lem:cv2}, the fact that $g_{(\xi,i,d_0)}$  have disjoint supports for  $\xi\neq\xi'$, and the fact that $\chi\tilde{\chi}=\chi$ we obtain 
\begin{align}
w_{q+1}^{(p,i)}\theta_{q+1}^{(p,i)}&=\sum_{n\geq3}\chi( \zeta|M^{i}_l|-n)\frac{n}{ \zeta}\sum_{\xi\in\Lambda^{n}}\Gamma_{\xi}\(\frac{M^{i}_l}{|M^{i}_l|}\)W_{(\xi,n,i)}\(\Theta_{(\xi,n,i)}-\mathbb{P}_{0}\Theta_{(\xi,n,i)}\)g_{(\xi,i,d_0)}g_{(\xi,i,d_0')}\notag\\ 
&=\sum_{n\geq3}\sum_{\xi\in\Lambda^{n}}\chi( \zeta|M^{i}_l|-n)\frac{n}{ \zeta}\Gamma_{\xi}\(\frac{M^{i}_l}{|M^{i}_l|}\)\mP_{\neq0}(W_{(\xi,n,i)}\Theta_{(\xi,n,i)})g_{(\xi,i,d_0)}g_{(\xi,i,d_0')}\notag\\ &\quad-\sum_{n\geq3}\chi( \zeta|M^{i}_l|-n)\frac{n}{ \zeta}\sum_{\xi\in\Lambda^{n}}\Gamma_{\xi}\(\frac{M^{i}_l}{|M^{i}_l|}\)W_{(\xi,n,i)}\mathbb{P}_{0}\Theta_{(\xi,n,i)}g_{(\xi,i,d_0)}g_{(\xi,i,d_0')}\notag\\
     &\quad +\sum_{n\geq3}\sum_{\xi\in\Lambda^{n}}\chi( \zeta|M_l^i|-n)\frac{n}{ \zeta}\Gamma_{\xi}\(\frac{M^{i}_l}{|M^{i}_l|}\)\xi (g_{(\xi,i,d_0)}g_{(\xi,i,d_0')}-1)\notag\\ 
     &\quad+\sum_{n\geq3}\chi( \zeta|M^{i}_l|-n)\frac{n}{ \zeta}\frac{M^{i}_l}{|M^{i}_l|}.\label{4wq+1ptheq+1p}
\end{align}

Now to deal with the undesired term in the last second line in \eqref{4wq+1ptheq+1p}, we define the oscillation perturbation as
\begin{align*}
 \theta_{q+1}^{(o,i)}:&=-\sigma^{-1}\sum_{n\geq3}\sum_{\xi\in\Lambda^{n}}h_{(\xi,i,d_0)}\div  \(\chi( \zeta|M_l^i|-n)\frac{n}{ \zeta}\Gamma_{\xi}\(\frac{M^{i}_l}{|M^{i}_l|}\)\xi\).
 \end{align*}
By  \eqref{eq:parth}, we have
\begin{align}
  \partial_t\theta^{(o,i)}_{q+1}&+ \sum_{n\geq3}\sum_{\xi\in\Lambda^{n}}(g_{(\xi,i,d_0)}g_{(\xi,i,d_0')}-1)\div\(\chi( \zeta|M_l^i|-n)\frac{n}{ \zeta}\Gamma_{\xi}\(\frac{M^{i}_l}{|M^{i}_l|}\)\xi\)\notag\\ 
  &=-\sigma^{-1}\sum_{n\geq3}\sum_{\xi\in\Lambda^{n}}h_{(\xi,i,d_0)}\partial_t\div \(\chi( \zeta|M_l^i|-n)\frac{n}{ \zeta}\Gamma_{\xi}\(\frac{M^{i}_l}{|M^{i}_l|}\)\xi\).\label{eq:4parttheo}
\end{align}
Moreover, since $M_l^i$ is compact supported and bounded, we know that $ \theta_{q+1}^{(o,i)}$ is mean-zero.

Furthermore,  we define the global perturbation as a solution to 
\begin{align}
     \partial_t \theta_{q+1}^{(g,i)}-\Delta\theta_{q+1}^{(g,i)}&=F_{q+1}^i,\notag\\ 
     \theta_{q+1}^{(g,i)}(0)&=0, \label{def:theq+1g}
\end{align}
where the force term $F^i_{q+1}$ is determined by the "inverse  divergence" of  stress term at step $q+1$ with support outside $\Omega_{q+1}$, see Proposition \ref{prop:fq+1} below for the precise definition. In order not to cause any confusion, we should remark that $F^i_{q+1}$ in fact depends only on $w_{q+1}^{(p,i)}$ and $\theta_{q+1}^{(p,i)}$. By  Proposition \ref{prop:fq+1}, $F^i_{q+1}$ is mean-zero, and satisfies $F_{q+1}^i(t)=0$ on $  [0,T_{q+1}]$. Then we have that $\theta^{(g,i)}_{q+1}$ is mean-zero, and satisfies $\theta_{q+1}^{(g,i)}(t)=0$ on $  [0,T_{q+1}]$.

Finally, the total perturbation and new iteration are defined by
\begin{align}
w_{q+1}:=\sum_{i=1}^{N_0}\(w_{q+1}^{(p,i)}+w_{q+1}^{(c,i)}\),\ \ v_{q+1}:=v_q+w_{q+1}.\notag
\end{align}
Moreover, 
for every $1\leq i\leq N_0$, we define
\begin{align}
\theta_{q+1}^{(loc,i)}:=\theta_{q+1}^{(p,i)}+\theta_{q+1}^{(c,i)}+\theta_{q+1}^{(o,i)},\ \ \theta^{i}_{q+1}:=\theta_{q+1}^{(loc,i)}+\theta_{q+1}^{(g,i)},
\ \ 
\rho^{i}_{q+1}:=\rho^{i}_q+\theta^{i}_{q+1}.\notag
\end{align}
 Then  $v_{q+1}$ is  mean-zero and divergence-free. By the definition, it is easy to see that $\int_{\mR^d}(\theta_{q+1}^{(p,i)}+\theta_{q+1}^{(c,i)})\dif x=0$, which implies that $\int_{\mR^d}\rho^{i}_{q+1}\dif x=0$.

Since $M^{i}_l(t)=0$ on $  [0,T_{q+1}]$,  we have $w_{q+1}(t)=\theta^i_{q+1}(t)=0$ on $[0,T_{q+1}]$. Then we have $v_{q+1}(t)=\rho_{q+1}^i(t)=0$ on $ [0,T_{q+1}]$. In summary we imply \eqref{bd:4mq=rhoq=0}   for $v_{q+1}$ and $\rho_{q+1}^i$. Similarly, since $\supp M_l^i\subset \Omega_{q+1}$, we have that $\supp w_{q+1}, \supp \theta_{q+1}^{(loc,i)}\subset \Omega_{q+1}$, which implies \eqref{bd:supp} for $v_{q+1}$.

Moreover, by the fact that $g_{(\xi,i,d_0)}$ have disjoint supports for distinct $i$, we have
\begin{align}
    w_{q+1}\theta^{i}_{q+1}=(w_{q+1}^{(p,i)}+w_{q+1}^{(c,i)})\theta^{(p,i)}_{q+1}+w_{q+1}(\theta^{(c,i)}_{q+1}+\theta^{(o,i)}_{q+1}+\theta_{q+1}^{(g,i)}).\label{bd:4wq+1rhoq+1i}
\end{align}

\subsection{Construction of the stress terms $M^{i}_{q+1}$}\label{sec:4defmq+1}
In this section, we present the exact expression for the stress term $M^{i}_{q+1}$. Unlike the torus setting, our construction cannot directly apply inverse divergence operators, since even when the stress term 
$M_l^i$   has compact support, its inverse divergence may  still remain support in the whole space. 

 \subsubsection{Inverse divergence iteration}\label{tamr}
First, we introduce the following inverse divergence iteration. A similar procedure for vector-valued functions was introduced in   \cite[Proposition A.17]{BMNV21}.

Let $\{\rho_{(n)}\}_{0\leq n\leq N+1}$ be the zero-mean smooth $\mathbb{T}^d$-periodic functions such that $\rho_{(n)} = \Delta\rho_{(n+1)}$. Then for any given function $G$ on $\R^d$,
 we have 
 $$G\rho_{(0)}=\div M_{(0)}+F_{(0)},$$
 where $M_{(0)}^i=G\partial_i\rho_{(1)},\ F_{(0)}:=-\sum_i\partial_iG\partial_i\rho_{(1)}.$

 Then, applying the decomposition to the error term $F_{(0)}$ at each step for $N$ times, we have 
 \begin{align}
     G\rho_{(0)}= \div M_{(N)}+ F_{(N)},
 \end{align}
 where for $1\leq i\leq d$
 $$M^{i}_{(N)} =\sum_{m=0}^N\sum_{|\alpha^i_m|=m+1}M^{(m+1)}_{\alpha^i_m}(G)\partial_{\alpha^i_m}\rho_{(m+1)},\ \ F_{(N)}=\sum_{|\alpha_N|=N+1}F^{(N+1)}_{\alpha_N}(G)\partial_{\alpha_N}\rho_{(N+1)}.$$
 Here $M^{(m+1)}_{\alpha_m}(G)$ is composed of the $m$-th derivative of $G$ and $F^{(m)}_{\alpha_m}(G)$
 the $m$-th derivative of $G$.
 
 Particularly, if $\supp G
 \subset [-\frac12,\frac12]^d$, we have $\supp M_{(N)},\supp F_{(N)}\subset\supp G$.

\subsubsection{Decomposition of $M^{i}_{q+1}$}
From now on, we will write $\Gamma_{\xi}(\frac{M^{i}_l}{|M^{i}_l|})=\Gamma_{\xi}$ if there  is no ambiguity. 
From   \eqref{bd:4wq+1rhoq+1i}, the definition of the perturbations
we obtain

\begin{align*}
 -\div M^{i}_{q+1}&=\partial_t\theta_{q+1}^i+\div (w_{q+1}^{(p,i)}\theta_{q+1}^{(p,i)}-M^{i}_l)-\div (M^{i}_q-M^{i}_l)\\
 & -\div ( \nabla\theta^{i}_{q+1})+\div(v_q\theta^{i}_{q+1}+w_{q+1}(\rho^{i}_q+\theta_{q+1}^{(c,i)}+\theta_{q+1}^{(o,i)}+\theta_{q+1}^{(g,i)}+\rho^{in,i})+w_{q+1}^{(c,i)}\theta_{q+1}^{(p,i)}),
\end{align*}
where together with \eqref{def:theq+1g}  we define  the linear error and commutator error by
\begin{align*} 
M^{i}_{lin}:&=-  \nabla\theta^{(loc,i)}_{q+1}+v_q\theta^{i}_{q+1}+w_{q+1}(\rho^{i}_q+\theta_{q+1}^{(c,i)}+\theta_{q+1}^{(o,i)}+\theta_{q+1}^{(g,i)}+\rho^{in,i}) +w_{q+1}^{(c,i)}\theta_{q+1}^{(p,i)},\\
M^{i}_{com}:&=M^{i}_q-M^{i}_l.
\end{align*}
Then it holds that
\begin{align*}
 -\div M^{i}_{q+1}&=\partial_t\theta_{q+1}^{(loc,i)}+\div (w_{q+1}^{(p,i)}\theta_{q+1}^{(p,i)}-M^{i}_l)+F_{q+1}^i+\div M^{i}_{lin}-\div M_{com}^i.
\end{align*}

For simplicity, we define
$\mathbb{P}_{\neq 0,\hat{\rho}^i_q}f:=f- \hat{\rho}^i_q\cdot\int_{\mathbb{R}^d}f\dif x,$
where $\hat{\rho}^i_q$ is defined in \eqref{def:thetaq+1c}.  
In particular, for any mean-zero function $f$, we have 
 $\mathbb{P}_{\neq 0,\hat{\rho}^i_q}f=f$. 
To define the oscillation error, by the identities \eqref{4wq+1ptheq+1p} and \eqref{eq:4parttheo} we have
\begin{align*}
 &\partial_t\theta^{(loc,i)}_{q+1}+\div (w_{q+1}^{(p,i)}\theta_{q+1}^{(p,i)}-M^{i}_l)\\
 &=\mathbb{P}_{\neq 0,\hat{\rho}^i_q} (\partial_t\theta^{(p,i)}_{q+1}+\div (w_{q+1}^{(p,i)}\theta_{q+1}^{(p,i)}-M^{i}_l)+\partial_t\theta^{(o,i)}_{q+1}
 )\\
&=\sum_{n\geq3}\sum_{\xi\in\Lambda^{n}}\mathbb{P}_{\neq0,\hat{\rho}^i_q}\(\partial_t[\chi( \zeta|M^{i}_l|-n)\(\frac{n}{ \zeta}\)^{1/d_0'}\Gamma_{\xi} 
g_{(\xi,i,d_0')}]\mathbb{P}_{\neq0}\Theta_{(\xi,n,i)}\)(:= \div M^{i}_{osc,t}+F^i_{osc,t})\notag\\ 
 &+\mathbb{P}_{\neq0,\hat{\rho}^i_q}\(\nabla[\chi( \zeta|M^{i}_l|-n)\frac{n}{ \zeta}\Gamma_{\xi} 
 ]g_{(\xi,i,d_0)}g_{(\xi,i,d_0')}\mP_{\neq0}(W_{(\xi,n,i)}\Theta_{(\xi,n,i)})\)(:= \div M^{i}_{osc,x}+F^i_{osc,x})\\
 &+\mathbb{P}_{\neq0,\hat{\rho}^i_q}\(\chi( \zeta|M^{i}_l|-n)\(\frac{n}{ \zeta}\)^{1/d_0'}g_{(\xi,i,d_0')}\Gamma_{\xi} \notag\\ &\quad\quad\times
 \(\partial_t\Theta_{(\xi,n,i)}+\(\frac{n}{ \zeta}\)^{1/d_0}g_{(\xi,i,d_0)}\div(W_{(\xi,n,i)}\Theta_{(\xi,n,i)})\)\)\\
 &-\div \(\sum_{n\geq3}\chi( \zeta|M^{i}_l|-n)\frac{n}{ \zeta}\sum_{\xi\in\Lambda^{n}}\Gamma_{\xi} 
 W_{(\xi,n,i)}\mathbb{P}_{0}\Theta_{(\xi,n,i)}g_{(\xi,i,d_0)}g_{(\xi,i,d_0')}\)(:= \div M^{i}_{osc,a})\notag\\
 &+\div \(\sum_{n\geq3}\chi( \zeta|M^{i}_l|-n)\frac{n}{ \zeta}\frac{M^{i}_l}{|M^{i}_l|}-M^{i}_l\)(:= \div M^{i}_{osc,c})\notag\\
 &-\div\(\sigma^{-1}\sum_{n\geq3}\sum_{\xi\in\Lambda^{n}}h_{(\xi,i,d_0)}\partial_t \(\chi(\zeta|M_l^i|-n)\frac{n}{ \zeta}\Gamma_{\xi}  
 \xi\)\)
 (:= \div M^{i}_{osc,o}),
\end{align*}
where the last forth term equals to 0 by \eqref{eq:4ptthe+n}.   Here we define 
\begin{align*}
M^{i}_{osc,a}&:=-\sum_{n\geq3}\chi( \zeta|M^{i}_l|-n)\frac{n}{ \zeta}\sum_{\xi\in\Lambda^{n}}\Gamma_{\xi}W_{(\xi,n,i)}\mathbb{P}_{0}\Theta_{(\xi,n,i)}g_{(\xi,i,d_0)}g_{(\xi,i,d_0')},\\
   M^{i}_{osc,c}&:=\sum_{n\geq3}\chi( \zeta|M^{i}_l|-n)\frac{n}{ \zeta}\frac{M^{i}_l}{|M^{i}_l|}- M^{i}_l,\\
    M^{i}_{osc,o}&:=-\sigma^{-1}\sum_{n\geq3}\sum_{\xi\in\Lambda^{n}}h_{(\xi,i,d_0)}\partial_t\(\chi(\zeta|M_l^i|-n)\frac{n}{ \zeta}\Gamma_{\xi}\xi\).
\end{align*}
 
Now, we  apply the inverse divergence iteration introduced above to define $M_{osc,x}$ and $M_{osc,t}$. We first write
\begin{align*}
    A^{1,(\xi,n,i)}_{q+1}&:= \partial_t[\chi( \zeta|M^{i}_l|-n)(\frac{n}{ \zeta})^{1/d_0'}\Gamma_{\xi} 
    g_{(\xi,i,d_0')}],\\ A^{2,(\xi,n,i)}_{q+1}&:=\nabla[\chi( \zeta|M^{i}_l|-n)\frac{n}{ \zeta}\Gamma_{\xi} 
    ]g_{(\xi,i,d_0)}g_{(\xi,i,d_0')}.
\end{align*}
 Then, we obverse that $\Theta_{(\xi,n,i)}$ and $W_{(\xi,n,i)}\Theta_{(\xi,n,i)}$  are both $(\mathbb{T}/r_\perp\lambda_{q+1})^d$-periodic, so it is possible to define $\Delta^{-n}\mathbb{P}_{\neq0}\Theta_{(\xi,n,i)}$ and $\Delta^{-n}\mP_{\neq0}(W_{(\xi,n,i)}\Theta_{(\xi,n,i)})$ for every $n\in\mN$.
By   using the inverse divergence iterations  in Section \ref{tamr} to $\rho_{(n)}=\Delta^{-n}\mathbb{P}_{\neq0}\Theta_{(\xi,n,i)}$ and to $\rho_{(n)}=\Delta^{-n}\mP_{\neq0}(W_{(\xi,n,i)}\Theta_{(\xi,n,i)})$ respectively we  obtain 
\begin{align*}
     A^{1,(\xi,n,i)}_{q+1}\mathbb{P}_{\neq0}\Theta_{(\xi,n,i)}=\div M^{1,(\xi,n,i)}_{(N)}+ F^{1,(\xi,n,i)}_{(N)},\\ A^{2,(\xi,n,i)}_{q+1}\mP_{\neq0}(W_{(\xi,n,i)}\Theta_{(\xi,n,i)})=\div M^{2,(\xi,n,i)}_{(N)}+ F^{2,(\xi,n,i)}_{(N)}.
\end{align*}
Here we recall that $N>4d$ is determined in Section \ref{sec:4choicepara}, and we abuse the notation to define
\begin{align*}
    (M^{1,(\xi,n,i)}_{(N)})^{j} &:=\sum_{m=0}^{N^2-1}\sum_{|\alpha_m^j|=m+1}M^{(m+1)}_{\alpha_m^j}(A_{q+1}^{1,(\xi,n,i)})\partial_{\alpha_m^j}\Delta^{-m-1}\mathbb{P}_{\neq0}\Theta_{(\xi,n,i)},\\
F^{1,(\xi,n,i)}_{(N)}&:=\sum_{|\alpha_N|=N^2}F^{(N^2)}_{\alpha_N}(A_{q+1}^{1,(\xi,n,i)})\partial_{\alpha_N}\Delta^{-N^2}\mathbb{P}_{\neq0}\Theta_{(\xi,n,i)},\\
(M^{2,(\xi,n,i)}_{(N)})^{j} &:=\sum_{m=0}^{N^2-1}\sum_{|\alpha_m^j|=m+1}M^{(m+1)}_{\alpha_m^j}(A_{q+1}^{2,(\xi,n,i)})\partial_{\alpha_m^j}\Delta^{-m-1}\mP_{\neq0}(W_{(\xi,n,i)}\Theta_{(\xi,n,i)}),\\
F^{2,(\xi,n,i)}_{(N)}&:=\sum_{|\alpha_N|=N^2}F^{(N^2)}_{\alpha_N}(A_{q+1}^{2,(\xi,n,i)})\partial_{\alpha_N}\Delta^{-N^2}\mP_{\neq0}(W_{(\xi,n,i)}\Theta_{(\xi,n,i)}).
\end{align*} Moreover, we have $\supp M^{1,(\xi,n,i)}_{(N)},\supp M^{2,(\xi,n,i)}_{(N)},\supp F^{1,(\xi,n,i)}_{(N)},\supp F^{2,(\xi,n,i)}_{(N)}\subset \Omega_{q+1}$, and  vanish on $[0,T_{q+1}]$.
Then, we could define  $$M^{i}_{osc,t}:=\sum_{n\geq3}\sum_{\xi\in\Lambda^{n}} M^{1,(\xi,n,i)}_{(N)},\ \ M^{i}_{osc,x}:=\sum_{n\geq3}\sum_{\xi\in\Lambda^{n}} M^{2,(\xi,n,i)}_{(N)},$$
and 
$$F^{i}_{osc,t}:=\sum_{n\geq3}\sum_{\xi\in\Lambda^{n}}\mathbb{P}_{\neq0,\hat{\rho}^i_q}F^{1,(\xi,n,i)}_{(N)},\ \ F^{i}_{osc,x}:=\sum_{n\geq3}\sum_{\xi\in\Lambda^{n}}\mathbb{P}_{\neq0,\hat{\rho}^i_q}F^{2,(\xi,n,i)}_{(N)}.$$

Finally, we define $M^{i}_{osc}:=M^{i}_{osc,t}+M^{i}_{osc,x}+M^{i}_{osc,a}+M^{i}_{osc,c}+M^{i}_{osc.o}$ and  
\begin{align}
  -{M}^i_{q+1}:=M^{i}_{osc}+M^{i}_{lin}-M_{com}^i,\ \   F^i_{q+1}:=-F^{i}_{osc,t}-F^{i}_{osc,x}.\label{def:Fq+1}
\end{align}
 
Since $M^{i}_l(t)=w_{q+1}(t)=\theta^i_{q+1}(t)=0$ on $[0,T_{q+1}]$,  we have $M_{q+1}^i(t)=0$ on $[0,T_{q+1}]$, which implies  \eqref{bd:4mq=rhoq=0} for $M_{q+1}^i$. Similarly, we have $\supp M_{q+1}^i\subset \Omega_{q+1}$, which implies \eqref{bd:supp} for $M_{q+1}^i$.

\subsection{Estimates of $w_{q+1}$}\label{sec:4estwq+1}
In this section, we aim to establish the desired estimates for the perturbation $w_{q+1}$. First, we establish the estimate of the amplitude functions which could be obtained by the same calculation as \cite[Proposition 5.2]{LRZ25}.
\bl\label{lem:4chi}
 For $ N\in\mathbb{N}_0,1\leq i\leq N_0$  we have
\begin{align*}
    \sum_{n\geq3}\|\chi( \zeta|M^{i}_l|-n)\|_{C_{t,x}^N}+ \sum_{n\geq3}\|\tilde{\chi}( \zeta|M^{i}_l|-n)\|_{C_{t,x}^N}&\lesssim l^{- (d+4)N-(d+2)},\\
     \sum_{n\geq3}\sum_{\xi\in\Lambda^{n}}\norm{\chi( \zeta|M^{i}_l|-n)\Gamma_\xi\(\frac{M^{i}_l}{|M^{i}_l|}\)}_{C_{t,x}^N}& \lesssim l^{-(2d+8)N-(d+2)},\\
   \(\frac{n}{ \zeta}\)^N1_{ \{\tilde{\chi}( \zeta|M^{i}_l|-n)>0\}}+\(\frac{n}{ \zeta}\)^N1_{ \{\chi( \zeta|M^{i}_l|-n)>0\}}&\lesssim l^{-N(d+2)},\\
     \sum_{n\geq3}\norm{\chi( \zeta|M^{i}_l|-n)\frac{M^{i}_l}{|M^{i}_l|}}_{C_{t,x}^N}&\lesssim l^{-(d+5)N-(2d+5)}.
\end{align*}
\el

 Then we introduce the improved H\"older’s inequality  by using the additional decorrelation between frequencies.
\bl$($\cite[Lemma A.4]{MNY24a}$)$\label{ihiot}
Let $d\geq2,p \in [1, \infty]$. Let $f : \mathbb{T}^d \to \mathbb{R}$ be a smooth function and $a$ be a smooth function on $\mR^d$ such that $\supp a\subset [-\frac12,\frac12]^d$. Then for any $\sigma\in\mathbb{N}$,
$$| \|a f(\sigma\cdot)\|_{L^p([-\frac12,\frac12]^d)}-\|a\|_{L^p([-\frac12,\frac12]^d)}\|f\|_{L^p(\mT^d)} |\lesssim\sigma^{-1/p}\|a\|_{C^1([-\frac12,\frac12]^d)}\|f\|_{L^p(\mT^d)}.$$
\el

 Recall from Section \ref{sec:4defq+1} that the perturbation $w_{q+1}$ decomposes as   $w_{q+1}=\sum_{i=1}^{N_0}(w_{q+1}^{(p,i)}+w_{q+1}^{(c,i)})$.  We first estimate the principle perturbations ${w}^{(p,i)}_{q+1}$ for $1\leq i\leq N_0$ in $L_t^{d_0}L^{d_0}$-norm. By Cauchy's inequality we have
\begin{align}
|w_{q+1}^{(p,i)}|^{d_0}&\leq\(\sum_{n\geq3}\tilde{\chi}( \zeta|M^{i}_l|-n)\)^{d_0-1}\sum_{n\geq3}\tilde{\chi}( \zeta|M^{i}_l|-n)\frac{n}{ \zeta}\Bigg{|}\sum_{\xi\in\Lambda^{n}}W_{(\xi,n,i)}g_{(\xi,i,d_0)}\bigg{|}^{d_0}\notag\\ 
&\lesssim \sum_{n\geq3}\tilde{\chi}( \zeta|M^{i}_l|-n)\frac{n}{ \zeta}\sum_{\xi\in\Lambda^{n}}\left|W_{(\xi,n,i)}g_{(\xi,i,d_0)}\right|^{d_0},\notag
\end{align}
where we used   the fact that $\sum_{n\in\mathbb{Z}} \tilde{\chi}(t - n) \leq2$. 

By applying the generalized H\"older inequality of 
Theorem \ref{ihiot} in spatial direction, together with the estimates for the building blocks in \eqref{int4} and Lemma \ref{lem:4chi} we deduce 
\begin{align}
\|w_{q+1}^{(p,i)}(t)\|^{d_0}_{L^{d_0}}
&\lesssim\sum_{n\geq3}\norm{\tilde{\chi}( \zeta|M^{i}_l(t)|-n)\frac{n}{ \zeta}}_{L^1}\sum_{\xi\in\Lambda^{n}}\|W_{(\xi,n,i)}\|_{C_tL^{d_0}(\mT^d)}^{d_0}g_{(\xi,i,d_0)}^{d_0}(t)\notag\\ 
&\quad+(r_\perp\lambda_{q+1})^{-1}\norm{\tilde{\chi}( \zeta|M^{i}_l(t)|-n)\frac{n}{ \zeta}}_{C_{t,x}^1}\sum_{\xi\in\Lambda^{n}}\|W_{(\xi,n,i)}\|_{C_tL^{d_0}(\mT^d)}^{d_0}g_{(\xi,i,d_0)}^{d_0}(t)\notag\\
&\lesssim\left(\norm{\sum_{n\geq3}\tilde{\chi}( \zeta|M^{i}_l(t)|-n)(M^{i}_l(t)+ \zeta^{-1}
)}_{L^1}+l^{-3d-8}\lambda_{q+1}^{-\frac{1}{N}}\right)\sum_{\xi\in\Lambda}g_{(\xi,i,d_0)}^{d_0}(t)\notag\\
&\lesssim (\|M^{i}_l(t)\|_{L^1}+ \delta_{q+1})\sum_{\xi\in\Lambda}g_{(\xi,i,d_0)}^{d_0}(t)
,\notag
\end{align}
where we used   the fact that $\sum_{n\in\mathbb{Z}} \tilde{\chi}(t - n) \leq2$, and used conditions on the parameters to have $(6d+1)\alpha-\frac{1}{N}<-\alpha<-2\beta$. Here we recall the notation $\Lambda=\Lambda^{1}\cup\Lambda^{2}$.
Then we apply the generalized H\"older inequality of Theorem \ref{ihiot} in time direction, the bounds \eqref{bd:gwnp} and \eqref{bd:4mlcn} to deduce  for some ${C}_v\geq1$
\begin{align}
\|w_{q+1}^{(p,i)}\|_{L^{d_0}_tL^{d_0}}^{d_0}&\lesssim (\|M^{i}_l\|_{L^1_tL^1}+ \delta_{q+1}+\sigma^{-1}\|M^{i}_l\|_{C_{t,x}^1})\sum_{\xi\in\Lambda}\norm{g_{(\xi,i,d_0)}}^{d_0}_{L_t^{d_0}}
\notag\\
    &\lesssim C_0 (\delta_{q+1}+\lambda_{q+1}^{(2d+6)\alpha-\frac{1}{2N}})\leq (\frac1{4N_0} C_v)^{d_0}C_0 \delta_{q+1},\label{bd:4wq+1pl2}
\end{align}
where we used  conditions on the parameters to have $(2d+6)\alpha-\frac{1}{2N}<-\alpha<-2\beta$.

For the general $L^u_tL^m$-norm with $u,m\in[1,\infty]$, by the estimates for the building blocks in \eqref{int2}-\eqref{int4}, \eqref{bd:gwnp}
and the estimate for the amplitude function in  Lemma \ref{lem:4chi} we obtain
\begin{align}
\|w_{q+1}^{(p,i)}\|_{L^u_tL^m}&\lesssim\sum_{n\geq3}\sum_{\xi\in\Lambda^{n}}\norm{\tilde{\chi}( \zeta|M^{i}_l|-n)\(\frac{n}{ \zeta}\)^{1/d_0}}_{C_{t,x}^0}\|W_{(\xi,n,i)}\|_{C_tL^m(\mT^d)}\|g_{(\xi,i,d_0)}\|_{L_t^u}\notag\\
&\lesssim l^{-2d-4}r_\perp^{\frac{d-1}{m}-\frac{d-1}{d_0}} r_\parallel^{\frac{1}{m}-\frac{1}{d_0}}\eta^{\frac1u-\frac1{d_0}},\label{bd:4wq+1plp}\\
\|w_{q+1}^{(c,i)}\|_{L^u_tL^m}&\lesssim\sum_{n\geq3}\sum_{\xi\in\Lambda^{n}}\norm{\tilde{\chi}( \zeta|M^{i}_l|-n)\(\frac{n}{ \zeta}\)^{1/d_0}}_{C^1_{t,x}}\notag\\ 
&\quad\quad \times\(\frac{1}{\lambda_{q+1}^2}\|\nabla\Phi_{(\xi,n,i)} \xi\cdot\nabla\psi_{(\xi,n,i)} \|_{C_tL^m(\mT^d)}+\|V_{(\xi,n,i)}\|_{C_tL^m(\mT^d)}\)\|g_{(\xi,i,d_0)}\|_{L_t^u}\notag\\
&
\lesssim l^{-3d-8}r_\perp^{\frac{d-1}{m}-\frac{d-1}{d_0}} r_\parallel^{\frac{1}{m}-\frac{1}{d_0}}\frac{r_\perp}{ r_\parallel}\eta^{\frac1u-\frac1{d_0}}.\label{bd:4wq+1clp}
\end{align}
Here and in the following, we should remark that since $\supp M_{l}^i\subset \Omega_{q+1}$, all integrations involving the building blocks reduce to computations on $[-\frac12,\frac12]^d$, equivalent to integration on the torus.

With these estimates, combining with the choice of parameters in  \eqref{para42} and the bound \eqref{bd:4wq+1pl2} we obtain 
\begin{align}
\|w_{q+1}\|_{L^{d_0}_tL^{d_0}}&\leq \frac{C_v}4C_0^{1/d_0}\delta_{q+1}^{1/d_0}+N_0Cl^{-3d-8}\frac{r_\perp}{ r_\parallel}\leq \frac{C_v}4C_0^{1/d_0}\delta_{q+1}^{1/d_0}+C\lambda_{q+1}^{(6d+17)\alpha-\frac1N}\notag \\
&\leq \frac{C_v}4C_0^{1/d_0}\delta_{q+1}^{1/d_0}+C\lambda_{q+1}^{-\alpha}\leq {C_v}C_0^{1/d_0}\delta_{q+1}^{1/d_0},\notag\end{align}
where we used  conditions on the parameters to have
${(6d+18)\alpha<\frac{1}{N}},N_0\lesssim\lambda_{q+1}^{\alpha}$ and chose $a$ large enough to absorb the universal constant. 
The above inequality 
yields
\eqref{bd:4vq+1-vqldd}  for $v_{q+1}$. Then \eqref{bd:4vql2} holds for $v_{q+1}$.

Then we turn to estimate the $L^r_tL^p$-norm and the $C_tL^s$-norm. 
Combining with the bounds \eqref{bd:4wq+1plp} and \eqref{bd:4wq+1clp} above we obtain 
\begin{align}
\|w_{q+1}\|_{L^{r}_tL^{p}}
\lesssim  N_0 l^{-3d-8}r_\perp^{\frac{d-1}{p}-\frac{d-1}{d_0}} r_\parallel^{\frac{1}{p}-\frac{1}{d_0}}\eta^{\frac1r-\frac1{d_0}}
\lesssim \lambda_{q+1}^{(6d+17)\alpha-\frac1N}
\lesssim\lambda_{q+1}^{-\alpha}
\leq\delta_{q+1}^{1/d_0},\label{bd:4wq+1lpr}\\
\|w_{q+1}\|_{C_tL^{s}}
\lesssim N_0  l^{-3d-8} r_\perp^{\frac{d-1}s-\frac{d-1}{d_0}} r_\parallel^{\frac1s-\frac{1}{d_0}}\eta^{-\frac1{d_0}}
\lesssim \lambda_{q+1}^{(6d+17)\alpha-\frac1N}\lesssim\lambda_{q+1}^{-\alpha}\leq\delta_{q+1}^{1/d_0},\notag
\end{align}
which leads to \eqref{bd:4vq+1-vqlpr}. Here we used \eqref{para42}, \eqref{bd:paralam4} and  
${(6d+18)\alpha<\frac{1}{N}},$ $N_0\lesssim\lambda_{q+1}^{\alpha}$. Then we chose $a$ large enough to absorb the universal constant. 

Next we estimate the $C_{t,x}^1$-norm. By the fact that 
\begin{align*}
    \partial_t(V_{(\xi,n,i)}(t))=\(\frac{n}{ \zeta}\)^{1/d_0}g_{(\xi,i,d_0)}\(\partial_t V_{(\xi)}\)\(\(\frac n{ \zeta}\)^{1/d_0}H_{(\xi,i,d_0)}(t)\),
\end{align*}
the estimates for the building blocks in \eqref{int4}, \eqref{bd:gwnp},  \eqref{eq:4wq+1p+wq+1c} and the estimates for the amplitude functions in Lemma \ref{lem:4chi} we have for $d_0\geq2$
\begin{align}
\|w_{q+1}\|_{C_{t,x}^1}&\lesssim\sum_{i=1}^{N_{0}}\sum_{n\geq3}\sum_{\xi\in\Lambda^{n}}\norm{\tilde{\chi}( \zeta|M^{i}_l|-n)\(\frac{n}{ \zeta}\)^{1/d_0}}_{C_{t,x}^{2}}\|\nabla V_{(\xi,n,i)}\|_{C_{t,x}^1}\|g_{(\xi,i,d_0)}\|_{C_t^1}\notag\\
&\lesssim \sum_{i=1}^{N_{0}}\sum_{n\geq3}\norm{\tilde{\chi}( \zeta|M^{i}_l|-n)\frac{n}{ \zeta}}_{C_{t,x}^{2}} \lambda_{q+1}\mu r_\parallel^{-\frac{1}{d_0}}r_\perp^{-\frac{d-1}{d_0}}\sigma\eta^{-1-\frac{2}{d_0}}\lesssim N_{0}\lambda_{q+1}^{(8d+24)\alpha+d+\frac72}.\notag
\end{align}
 Thus by $N_{0}\lesssim\lambda_{q+1}^{\alpha},(8d+25)\alpha<\frac{1}{2}$ we obtain
 \begin{align*}
     \|v_{q+1}\|_{C_{t,x}^1}\leq\|v_{q}\|_{C_{t,x}^1}+\|w_{q+1}\|_{C_{t,x}^1}\leq C_0^{1/d_0}\lambda_{q}^{d+4}+\frac12\lambda_{q+1}^{d+4}\leq C_0^{1/d_0}\lambda_{q+1}^{d+4},
 \end{align*}
which implies \eqref{bd:4rhoqc1} for $v_{q+1}$. Here we chose $a$ large enough to absorb the universal constant.

\subsection{Estimates of $\theta^{(loc,i)}_{q+1}$ }\label{sec:4esttheq+1}

Similarly as before, we first estimate the principle perturbations  $\theta^{(p,i)}_{q+1}$ in $L^{d_0'}_tL^{d_0'}$.  By  the fact that $\Gamma_{\xi}$ are uniformly bounded we have
\begin{align}
|\theta_{q+1}^{(p,i)}|^{d_0'}
&\lesssim \sum_{n\geq3}\chi( \zeta|M^{i}_l|-n)\frac{n}{ \zeta}\sum_{\xi\in\Lambda^{n}}\left|\mP_{\neq0}\Theta_{(\xi,n,i)}g_{(\xi,i,d_0')}\right|^{d_0'}.\notag
\end{align}
 Then by the same argument as in \eqref{bd:4wq+1pl2}, we have for some ${C}_\rho\geq1$
\begin{align}
    \|\theta_{q+1}^{(p,i)}\|_{L^{d_0'}_tL^{d_0'}}^{d_0'}&\lesssim \(\|M^{i}_l\|_{L^1_tL^1}+ \delta_{q+1}+\sigma^{-1}\|M^{i}_l\|_{C_{t,x}^1}\)
    \sum_{\xi\in\Lambda }\norm{g_{(\xi,i,d_0')}}^{d_0'}_{L_t^{d_0'}}
    \notag\\
    &\lesssim  C_0(\delta_{q+1}+\lambda_{q+1}^{(2d+6)\alpha-\frac{1}{2N}})\leq  (\frac12C_\rho)^{d_0'}C_0\delta_{q+1},\label{bd:4theq+1pl2}
\end{align}
where we used  conditions on the parameters to have $(2d+6)\alpha-\frac{1}{2N}<-\alpha<-2\beta$.

For general $L_t^uL^m$-norm with $u,m\in[1,\infty]$, by the estimates for the building blocks in \eqref{int4theta}, \eqref{bd:gwnp} and Lemma \ref{lem:4chi} we obtain
\begin{align}
\|\theta_{q+1}^{(p,i)} \|_{L^u_tL^{m}}&\lesssim\sum_{n\geq3}\sum_{\xi\in\Lambda^{n}}\norm{\chi( \zeta|M^{i}_l|-n)\(\frac{n}{ \zeta}\)^{1/d_0'}\Gamma_{\xi} 
}_{C_{t,x}^0}\|\Theta_{(\xi,n,i)} \|_{C_tL^{m}(\mT^d)}\|g_{(\xi,i,d_0')}\|_{L_t^{u}}\notag\\ 
&\lesssim l^{-2d-4}r_\perp^{\frac{d-1}{m}-\frac{d-1}{d_0'}} r_\parallel^{\frac{1}{m}-\frac{1}{d_0'}}\eta^{\frac1u-\frac{1}{d_0'}}.\label{bd:4theq+1plp}
\end{align}

Now we   estimate the $C_{t,x}^1$-norm. 
 By \eqref{int4theta}, \eqref{bd:gwnp}, \eqref{def:4para}, \eqref{bd:paralam4}  and Lemma \ref{lem:4chi} we have
\begin{align}
\|\theta_{q+1}^{(p,i)} \|_{C_{t,x}^1}
&\lesssim\sum_{n\geq3}\sum_{\xi\in\Lambda^{n}}\norm{\chi( \zeta|M^{i}_l|-n)\(\frac{n}{ \zeta}\)^{1/d_0'}\Gamma_{\xi} 
}_{C_{t,x}^1}\|\Theta_{(\xi,n,i)} \|_{C_{t,x}^1}\|g_{(\xi,i,d_0')}\|_{C_t^1}\notag\\ 
&\lesssim \sum_{n\geq3} \norm{\chi( \zeta|M^{i}_l|-n)\frac{n}{ \zeta}\Gamma_{\xi}}_{C_{t,x}^1}\lambda_{q+1}\mu r_\parallel^{-\frac{1}{d_0'}}r_\perp^{-\frac{d-1}{d_0'}}\sigma\eta^{-1-\frac{1}{d_0}-\frac{1}{d_0'}}\lesssim
\lambda_{q+1}^{(8d+24)\alpha+d+\frac{7}{2}},\notag\\
\|\theta_{q+1}^{(c,i)} \|_{C_{t,x}^1}&\lesssim\|\theta_{q+1}^{(p,i)} \|_{C_tL^1}
\lesssim l^{-2d-4}r_\perp^{d-1-\frac{d-1}{d_0'}} r_\parallel^{1-\frac{1}{d_0'}}\eta^{-\frac{1}{d_0'}}\lesssim  \lambda_{q+1}^{(4d+8)\alpha-\frac1N}\lesssim  \lambda_{q+1}^{-\alpha},\label{bd:4theq+1clp}\\
     \|   \theta_{q+1}^{(o,i)}\|_{C_{t,x}^1}&\lesssim \sigma^{-1}\sum_{n\geq3}\sum_{\xi\in\Lambda^{n}}\|h_{(\xi,i,d_0)}\|_{C_t^1}\norm{\div\(\chi( \zeta|M_l|-n)\frac{n}{ \zeta}\Gamma_{\xi}\xi\)}_{C_{t,x}^1}\notag\\ 
     &\lesssim \sigma^{-1}\eta^{-1} l^{-6d-20}\lesssim \lambda_{q+1}^{(12d+40)\alpha+2},\notag\\\
\|\theta_{q+1}^{(o,i)}\|_{C_tC^1}&\lesssim \sigma^{-1}\sum_{n\geq3}\sum_{\xi\in\Lambda^{n}}\|h_{(\xi,i,d_0)}\|_{L_t^\infty}\norm{\chi( \zeta|M_l|-n)\frac{n}{ \zeta}\Gamma_{\xi}  
\xi}_{C_tC^2}\notag\\ 
&\lesssim \sigma^{-1} l^{-6d-20}\lesssim \lambda_{q+1}^{(12d+40)\alpha-\frac1{2N}}\lesssim  \lambda_{q+1}^{-\alpha},\label{bd:4theq+1olp}
\end{align}
where we used  conditions on the parameters to have $(12d+41)\alpha<\frac1{2N}$ and chose $a$ large enough to absorb the universal constant.

Moreover,   once we  take a space derivative of $\theta_{q+1}^{(p,i)}$ or $\theta_{q+1}^{(o,i)}$, by  Lemma \ref{lem:4chi}, \eqref{int4theta} and \eqref{bd:gwnp}, we will obtain an extra power of $\lambda_{q+1}$, which implies that
\begin{align}    \|\theta_{q+1}^{(loc,i)} \|_{C_{t,x}^1}+  \lambda_{q+1}^{-1}\|\nabla\theta_{q+1}^{(loc,i)} \|_{C_{t,x}^1}\lesssim \lambda_{q+1}^{(12d+40)\alpha+d+\frac{7}{2}}
.\label{bs:theq+1locc2}\end{align}

To end this section, we estimate $\theta_{q+1}^{(p,i)}$ in $W^{1,1+\epsilon}$-norm for some $\epsilon>0$ small enough, which will be used below in the estimate of stress term. By the bounds for the building blocks in \eqref{int4theta}, \eqref{bd:gwnp}  and the bounds for the amplitude functions in Lemma \ref{lem:4chi} we have
\begin{align}
\|\theta_{q+1}^{(p,i)}\|_{L^1_tW^{1,1+\epsilon}}&\lesssim\sum_{n\geq3}\sum_{\xi\in\Lambda^{n}}\norm{\chi( \zeta|M^{i}_l|-n)\(\frac{n}{ \zeta}\)^{1/d_0'}\Gamma_{\xi}}_{C_{t,x}^1}\|\Theta_{(\xi,n,i)} \|_{C_tW^{1,1+\epsilon}(\mT^d)}\|g_{(\xi,i,d_0')}\|_{L_t^1}\notag\\ 
&\lesssim l^{-4d-12}\lambda_{q+1}r_\parallel^{\frac1{1+\epsilon}-\frac{1}{d_0'}}r_\perp^{\frac{d-1}{1+\epsilon}-\frac{d-1}{d_0'}}\eta^{1-\frac1{d_0'}}
\lesssim\lambda_{q+1}^{(8d+24)\alpha-\frac1N+d\epsilon}\lesssim \lambda_{q+1}^{-\alpha} ,\label{bd:4theq+1w1s}
\end{align}
where we used the choice of parameters in \eqref{bd:paralam4} and  chose $\epsilon>0$ small enough such that $d\epsilon<\alpha$. We also used  conditions on the parameters to have $(8d+26)\alpha<\frac1N$.

\subsection{Estimates of $F^{i}_{q+1}$}\label{sec:4estFq+1}
We recall the functions $F_{q+1}^i$ defined in \eqref{def:Fq+1}. 
In this section, we provide estimates for $F_{q+1}^i$, which will be used 
to derive the corresponding estimates for $\theta_{q+1}^{(g,i)}$.

\bp\label{prop:fq+1} The functions 
$F_{q+1}^i$ defined in \eqref{def:Fq+1} satisfy the following properties:   $F_{q+1}^i$ are mean-zero,
$\supp F_{q+1}^i\subset \Omega_{q+1}$ and $F_{q+1}^i=0 $ on $[0,T_{q+1}]$. Moreover, it holds that for $k,j\in\mN$ \begin{align}
  \|\partial_t^k\nabla^j F_{q+1}^i\|_{C^0_{t,x}}\lesssim \lambda_{q+1}^{2k+j-1}.\notag
\end{align}
\ep
\begin{proof}The first three properties follow directly from the definition. To establish the H\"older continuity, we begin by considering the case $k=j=0$
  and analyze each term separately. 
 We  obverse that $\Theta_{(\xi,n,i)}$ and $W_{(\xi,n,i)}\Theta_{(\xi,n,i)}$  are both $(\mathbb{T}/r_\perp\lambda_{q+1})^d$-periodic, and $\nabla^m\Delta^{-m}$ are  $(-m)$-homogeneous operators.  By the estimates for the amplitude functions and for the building blocks in  Lemma \ref{lem:4chi}, \eqref{int4theta}, \eqref{bd:gwnp} and \eqref{bd:paralam4} respectively, we have
 \begin{align*}
  \| F^{i}_{osc,t}\|_{C_{t,x}^0}&\lesssim   \sum_{n\geq3}\sum_{\xi\in\Lambda^{n}}\|F^{1,(\xi,n,i)}_{(N)}\|_{C_{t,x}^0}\notag\\ 
  &\lesssim \sum_{n\geq3}\sum_{\xi\in\Lambda^{n}} \|\nabla^{N^2}\partial_t[\chi( \zeta|M^{i}_l|-n)(\frac{n}{ \zeta})^{1/d_0'}\Gamma_{\xi} g_{(\xi,i,d_0')}]\|_{C_{t,x}^0}\|\nabla^{N^2}\Delta^{-N^2}\mathbb{P}_{\neq0}\Theta_{(\xi,n,i)}\|_{C_{t,x}^0}\\
  &\lesssim l^{-4d-12-(2d+8)N^2}(r_\perp^{d-1-\frac{d-1}{d_0'}}r_\parallel^{1-\frac{1}{d_0'}}\sigma \eta^{-\frac{1}{d_0'}})r_\perp^{-d+1}r_\parallel^{-1}\eta^{-1}(r_\perp\lambda_{q+1})^{-N^2}\\
  &\lesssim \lambda_{q+1}^{(8d+24+(4d+16)N^2)\alpha+d+1-N}\lesssim \lambda_{q+1}^{d+2-3N/4}\lesssim  \lambda_{q+1}^{-1},\\
\|F^{i}_{osc,x}\|_{C_{t,x}^0}&\lesssim \sum_{n\geq3}\sum_{\xi\in\Lambda^{n}}\|F^{2,(\xi,n,i)}_{(N)}\|_{C_{t,x}^0}\notag\\ 
&\lesssim \sum_{n\geq3}\sum_{\xi\in\Lambda^{n}}\|\nabla^{N^2+1}[\chi( \zeta|M^{i}_l|-n)\frac{n}{ \zeta}\Gamma_{\xi}]g_{(\xi,i,d_0)}g_{(\xi,i,d_0')}\|_{C_{t,x}^0}\\
&\quad\quad\quad\quad\quad\quad\times\|\nabla^{N^2}\Delta^{-N^2}\mP_{\neq0}(W_{(\xi,n,i)}\Theta_{(\xi,n,i)})\|_{C_{t,x}^0}\\
&\lesssim l^{-4d-12-(2d+8)N^2}\eta^{-1}r_\perp^{-d+1}r_\parallel^{-1}(r_\perp\lambda_{q+1})^{-N^2}\\
&\lesssim \lambda_{q+1}^{(8d+24+(4d+16)N^2)\alpha+d+1-N}\lesssim \lambda_{q+1}^{d+2-3N/4}
\lesssim  \lambda_{q+1}^{-1},
 \end{align*}
 where we used to conditions on the parameters to deduce that $(8d+32)\alpha<\frac1{2N}$.

When  taking a space derivative on $F_{q+1}^i$, by  Lemma \ref{lem:4chi}, \eqref{int4theta} and \eqref{bd:gwnp}, we will obtain an extra power of $l^{-2d-8}\lesssim\lambda_{q+1} $ from the amplitude functions and a extra power of $\lambda_{q+1}$ from the building blocks. Similarly, when taking a time derivative,   we will obtain an extra power of $l^{-2d-8}+\sigma\eta^{-1}+\frac{r_\perp\lambda\mu}{r_\parallel}\lesssim\lambda_{q+1}^2$.  Analogous to the previous analysis, we obtain the estimate on the derivatives.
\end{proof}

\subsection{Estimates of $\theta^{(g,i)}_{q+1}$}\label{sec:4esttheq+1g}
Since the functions 
$\theta_{q+1}^{(g,i)}$
  vanish identically on the interval 
$  [0,T_{q+1}]$, our subsequent analysis can be restricted to the regime $t\geq T_{q+1}\geq\frac16$.
We begin by establishing a key estimate for heat kernel convolutions:
\bp\label{prop:heatkernel}
Let $d\geq2$. For any $t\in[\frac1{12},1]$ and $x\in\mR^d$, we have
$$\int_0^t\int_{[-\frac12,\frac12]^d}p( s,x-y)\dif s\dif y\lesssim t\int_{[-\frac12,\frac12]^d}p(t,x-y)\dif y.$$
\ep

\begin{proof}
If $|x|\geq9$, we have $|x-y|^2/4t\geq 1$, then it
  suffices to prove that $\int_0^tp(s,x)\dif s\lesssim tp(t,x)$ for $a:=|x|^2/4t\geq 1$. By change of variable we have that \begin{align*}
 \int_0^t(4\pi s)^{-d/2}&e^{-|x|^2/4s}\dif s=\int_{|x|^2/4t}^\infty(\pi|x|^2/u)^{-d/2}e^{-u}\frac{|x|^2}{4u^2}\dif u\notag\\ 
    &=t(4\pi t)^{-d/2}a^{-d/2+1}\int_{a}^\infty u^{d/2-2}e^{-u}\dif u\leq t(4\pi t)^{-d/2}\int_{a}^\infty (u/a)^{d/2-1}e^{-u}\dif u\notag\\ 
    &\lesssim t(4\pi t)^{-d/2}\int_{a}^\infty (1+u-a)^{d/2-1}e^{-u}\dif u\leq t(4\pi t)^{-d/2}e^{-a}\int_{1}^\infty u^{d/2-1}e^{-u}\dif u\lesssim tp(t,x).
\end{align*}
Otherwise if $|x|<9 $, since $t\geq1/12$, we have 
\begin{align*}
    \int_{[-\frac12,\frac12]^d}p(t,x-y)\dif y\geq\int_{[-\frac12,\frac12]^d}(4\pi t)^{-d/2}e^{-25/t}\dif y\geq (4\pi)^{-d/2}e^{-300},
\end{align*}
which implies that
\begin{align*}
  \int_0^t\int_{[-\frac12,\frac12]^d}p(s,x-y)\dif s\dif y\leq t \lesssim t\int_{[-\frac12,\frac12]^d}p(t,x-y)\dif y.
\end{align*}

\end{proof}

Then, we are in position to bound $\theta_{q+1}^{(g,i)}$. By   Proposition \ref{prop:fq+1} and Proposition \ref{prop:heatkernel}, we have that for $t-\frac1{12}\geq \frac1{12}$

\begin{align}\label{bd:theq+1gic0}
   | \theta_{q+1}^{(g,i)}(t,x)|&\leq \int_0^t\int_{\mR^d}p(t-s,x-y)|F_{q+1}^i(s,y)|\dif s\dif y\lesssim  \lambda_{q+1}^{-1}\int_{0}^{t-\frac1{12}}\int_{[-\frac12,\frac12]^d}p(s,x-y)\dif s\dif y\notag\\
   &\lesssim  \lambda_{q+1}^{-1}\int_{[-\frac12,\frac12]^d}p(t-\frac1{12},x-y)\dif y\lesssim  \lambda_{q+1}^{-1}\int_{[-\frac12,\frac12]^d}e^{-c|x-y|^2}\dif y\lesssim  \lambda_{q+1}^{-1}e^{-c|x|^2},
\end{align} 
which implies that
\begin{align}
    \|\theta^{(g,i)}_{q+1}\|_{L^{d_0'}_tL^{d_0'}}+ \|\theta^{(g,i)}_{q+1}\|_{C_tL^1}+\|\theta^{(g,i)}_{q+1}\|_{C_{t,x}^0}\lesssim \lambda_{q+1}^{-1}.\label{bd:theq+1g}
\end{align}
Then together with the estimates on $\theta_{q+1}^{(loc,i)}$ in Section \ref{sec:4esttheq+1} and the choice of parameters in \eqref{para42}, \eqref{bd:paralam4}    we imply 
\begin{align}
\|\rho^{i}_{q+1}-\rho^{i}_q\|_{L^{d_0'}_tL^{d_0'}}
&\leq\|\theta^{i}_{q+1}\|_{L^{d_0'}_tL^{d_0'}}
\leq \frac{C_\rho}2C_0^{1/{d_0'}}\delta_{q+1}^{1/{d_0'}}+C\lambda_{q+1}^{-\alpha}\leq C_\rho C_0^{1/{d_0'}}\delta_{q+1}^{1/{d_0'}},\notag\\
\|\rho^{i}_{q+1}-\rho^{i}_q\|_{C_tL^1}
&\leq\|\theta^{i}_{q+1}\|_{C_tL^1}
\lesssim \lambda_{q+1}^{(4d+8)\alpha-\frac1{N}}+\lambda_{q+1}^{-\alpha}\lesssim  \lambda_{q+1}^{-\alpha}\leq  \delta_{q+1}^{1/d_0'},\label{bd:theq1l1}
\end{align}
which implies \eqref{bd:4vq+1-vqldd}, \eqref{bd:4rhoq+1-rhoql1}  and then \eqref{bd:4vql2} for $\rho^{i}_{q+1}$. Here we used $(4d+9)\alpha<\frac1{N}$, and chose  $a$ large enough to absorb the universal constant. 

Next, by the fact that $\supp\theta_{q+1}^{(c,i)},\supp\theta
_{q+1}^{(o,i)}\subset [\frac16,1]\times [-\frac12,\frac12]^d$, together with  \eqref{bd:4theq+1clp}, \eqref{bd:4theq+1olp} we know that 
\begin{align*}
  \theta_{q+1}^{(c,i)}(t,x)+\theta
_{q+1}^{(o,i)}(t,x)\geq - \lambda_{q+1}^{-\alpha}1_{\{[\frac16,1]\times [-\frac12,\frac12]^d\}}(t,x)\geq -\frac{1}{c_d}\lambda_{q+1}^{-\alpha}\overline p(t-\frac1{12},x).
\end{align*}
By the fact that $\theta_{q+1}^{(p,i)}$ is non-negative,  and the choice of parameters in \eqref{para42}, \eqref{bd:theq+1gic0}  we have
\begin{align}
    (\rho^{i}_{q+1} - \rho^{i}_q)(t,x)&\geq-
C \lambda_{q+1}^{-\alpha}\overline p(t-\frac1{12},x)\geq - \delta_{q+1}^{1/d_0'}\overline p(t-\frac1{12},x),\notag
\end{align}
which yields \eqref{bd:4rhoq+1-rhoqgeq}.  Here we chose $a$ large enough to absorb the universal constant. 

By a similar calculation as \eqref{bd:theq+1gic0} and Proposition \ref{prop:fq+1}, we obtain that
\begin{align}
 | \partial_t\theta_{q+1}^{(g,i)}|+  | \nabla\theta_{q+1}^{(g,i)}|&\leq \int_0^t\int_{\mR^d}p(t-s,x-y)(| \partial_t F_{q+1}^i(s,y)|+| \nabla F_{q+1}^i(s,y)|)\dif s\dif y\lesssim \lambda_{q+1}^2.\notag
\end{align}
Then together with  \eqref{bs:theq+1locc2}     we imply 
\begin{align}   \|\rho^{i}_{q+1}\|_{C_{t,x}^1}\leq\|\rho^{i}_{q}\|_{C_{t,x}^1}+\|\theta^{i}_{q+1}\|_{C_{t,x}^1}
    \leq C_0^{1/d_0'}\lambda_{q}^{d+4}+ \frac12\lambda_{q+1}^{d+4}\leq C_0^{1/d_0'} \lambda_{q+1}^{d+4},\notag
\end{align}
which implies  \eqref{bd:4rhoqc1} for $\rho^{i}_{q+1}$.
Here we chose $a$ large enough to absorb the universal constant.

\subsection{Estimates of $M^{i}_{q+1}$}\label{sec:4estmq+1}
We estimate each terms in the definition of $M^{i}_{q+1}$ separately. 
 \subsubsection{Estimate of oscillation error $M_{osc}^i$.}
 Since $\Theta_{(\xi,n,i)}$ and $W_{(\xi,n,i)}\Theta_{(\xi,n,i)}$  are both $(\mathbb{T}/r_\perp\lambda_{q+1})^d$-periodic, by the estimates for the amplitude functions and for the building blocks in  Lemma \ref{lem:4chi}, \eqref{int4theta} and \eqref{bd:gwnp} respectively   we obtain
\begin{align*}
  \| M^{i}_{osc,x}\|_{L^1_tL^1}&\lesssim \sum_{n\geq3}\sum_{\xi\in\Lambda^{n}} \sum_{m=0}^{N^2-1}\| \nabla
  ^m\partial_t[\chi( \zeta|M^{i}_l|-n)(\frac{n}{ \zeta})^{1/d_0'}\Gamma_{\xi} 
    g_{(\xi,i,d_0')}]\|_{L_t^1C^0}\\
    &\quad\quad\quad\times(r_\perp\lambda_{q+1})^{-m-1}\|\Theta_{(\xi,n,i)}\|_{C_tL^1(\mT^d)},\\
    &\lesssim l^{-4d-12}r_\perp^{d-1-\frac{d-1}{d_0'}}r_\parallel^{1-\frac{1}{d_0'}}\sigma \eta^{-\frac{1}{d_0'}}\sum_{m=0}^{N^2-1}\lambda_{q+1}^{((4d+16)\alpha-\frac1N)m-\frac1N}\lesssim \lambda_{q+1}^{(8d+24)\alpha-\frac 1{N}}\lesssim  \lambda_{q+1}^{-\alpha},\\
 \| M^{i}_{osc,t}\|_{L^1_tL^1}&\lesssim \sum_{n\geq3}\sum_{\xi\in\Lambda^{n}}\sum_{m=0}^{N^2-1}\|\nabla^{1+m}[\chi( \zeta|M^{i}_l|-n)\frac{n}{ \zeta}\Gamma_{\xi} 
    ]\|_{C_{t,x}^0}\|g_{(\xi,i,d_0)}g_{(\xi,i,d_0')}\|_{L_t^1}\\
    &\quad\quad\quad\times(r_\perp\lambda_{q+1})^{-m-1}\|W_{(\xi,n,i)}\Theta_{(\xi,n,i)}\|_{C_tL^1(\mT^d)},\\
    &\lesssim  l^{-4d-12}\lambda_{q+1}^{-\frac1N}\sum_{m=0}^{N^2-1}\lambda_{q+1}^{((4d+16)\alpha-\frac1N)m}\lesssim \lambda_{q+1}^{(8d+24)\alpha-\frac{1}{N}}
\lesssim  \lambda_{q+1}^{-\alpha},
\end{align*}
where we used   conditions on the parameters to have $(8d+25)\alpha<\frac{1}{N}$.

By Lemma \ref{lem:4chi}, \eqref{int4theta} and \eqref{bd:gwnp}  again we obtain
\begin{align*}
   \| M^{i}_{osc,a}\|_{L^1_tL^1}&\lesssim \sum_{n\geq3}\sum_{\xi\in\Lambda^{n}}\norm{\chi( \zeta|M^{i}_l|-n)\frac{n}{ \zeta}\Gamma_{\xi}}_{C_{t,x}^0}\|W_{(\xi,n,i)}\|_{C_tL^1(\mT^d)}\|\Theta_{(\xi,n,i)}\|_{C_tL^1(\mT^d)}\\
   &\lesssim l^{-2d-4}r_\perp^{d-1} r_\parallel\lesssim \lambda_{q+1}^{(4d+8)\alpha-\frac1N}\lesssim  \lambda_{q+1}^{-\alpha}, 
\end{align*}
where we used   conditions on the parameters to have $(4d+9)\alpha<\frac{1}{N}$.

The stress term $M^{i}_{osc,c}$ is bounded  similarly as in \cite[Section 5.4.2]{LRZ25}. 
\begin{align}
     \left| M_{osc,c}^i\right|&\leq \left|  \sum_{n=-1}^{2}\chi( \zeta|M_l^i|-n)M_l^i\right|+\left|\sum_{n\geq3}\chi( \zeta|M_l^i|-n)(\frac{n}{ \zeta}\frac{M_l^i}{|M_l^i|}- M_l^i)\right|\notag\\
     &\leq \frac{3}{ \zeta}+\sum_{n\geq3}\chi( \zeta|M_l^i|-n)\left|\frac{n}{ \zeta}- |M_l^i|\right|\leq \frac{3}{20} \delta_{q+2}+\frac{1}{20} \delta_{q+2}
     \leq \frac15C_0 \delta_{q+2}.\notag
\end{align}
 By the bounds  \eqref{bd:gwnp}, \eqref{bd:4rql1} and \eqref{para42}  we have
\begin{align}
   \|  M^{i}_{osc,o}\|_{L^1_tL^1}&\lesssim \sigma^{-1}\sum_{n\geq3}\sum_{\xi\in\Lambda^{n}}\|h_{(\xi,d_0)}\|_{L_t^\infty}\norm{\partial_t[\chi(\zeta|M_l^i|-n)\frac{n}{ \zeta}\Gamma_{\xi}  
   \xi]}_{C_{t,x}^0}\notag\\ 
   &\lesssim  \sigma^{-1}l^{-4d-12}
   \lesssim \lambda_{q+1}^{(8d+24)\alpha-\frac1{2N}}
   \lesssim  \lambda_{q+1}^{-\alpha},\notag
\end{align}
where we used   conditions on the parameters to have $(8d+25)\alpha<\frac{1}{2N}$. 

In summary, we have 
\begin{align}
        \|  M^{i}_{osc}\|_{L^1_tL^1}\leq CC_0\lambda_{q+1}^{-\alpha}+\frac15C_0 \delta_{q+2}\leq \frac13C_0 \delta_{q+2},\label{bd:moscl1}
\end{align}
where we choose $a$ large enough to absorb the constant.

Moreover, similar to the case of $F_{q+1}^i$,  by  Lemma \ref{lem:4chi}, \eqref{int4theta} and \eqref{bd:gwnp}, if we take a space/time derivative on $M_{osc}^i-M_{osc,c}^i$, we will obtain an extra power   of $\lambda_{q+1}^2$ at most, so we have
 \begin{align*}
    \| \partial_t(M_{osc}^i-M_{osc,c}^i)\|_{L^1_tL^1}+   \|  \nabla (M_{osc}^i-M_{osc,c}^i)\|_{L^1_tL^1}\lesssim  \lambda_{q+1}^{2-\alpha}.
 \end{align*}
 As for $M_{osc,c}^i$, by Lemma \ref{lem:4chi} we have 
 \begin{align*}
     \|M_{osc,c}^i\|_{C_{t,x}^1}\lesssim \|M_{q}^i\|_{C_{t,x}^1}+\sum_{n\geq3}\norm{\chi( \zeta|M^{i}_l|-n)\frac{n}{ \zeta}\frac{M^{i}_l}{|M^{i}_l|}}_{C_{t,x}^1} \lesssim \lambda_q^{2d+8}+ l^{-(4d+12)}\lesssim \lambda_{q+1}^{2d+7},
 \end{align*}
 which implies that 
 \begin{align}
    \| \partial_tM_{osc}^i\|_{L^1_tL^1}+   \|  \nabla M_{osc}^i\|_{L^1_tL^1}\lesssim  \lambda_{q+1}^{2d+7}\leq \frac13\lambda_{q+1}^{2d+8},\label{bd:moscc1}
 \end{align}
 where we chose $a$ large enough to absorb the universal constant.

 \subsubsection{Estimate of linear error $M^{i}_{lin}$.}
 By  the bounds in  \eqref{bd:4theq+1olp}-\eqref{bd:4theq+1w1s}   we have
\begin{align*}
     \| \nabla\theta^{(loc,i)}_{q+1}\|_{L^1_tL^1}
     &\lesssim \|\theta^{(p,i)}_{q+1}\|_{L^1_tW^{1+\epsilon}}+\|\theta^{(o,i)}_{q+1}\|_{L^1_tW^{1+\epsilon}} \lesssim   \lambda_{q+1}^{-\alpha}.
\end{align*}
By the estimates
in \eqref{bd:4rhoqc1}, \eqref{bd:4wq+1lpr}, \eqref{bd:4theq+1clp}, \eqref{bd:4theq+1olp},    \eqref{bd:theq+1g} and \eqref{bd:theq1l1} respectively we obtain
\begin{align}
\|v_q&\theta^{i}_{q+1}+w_{q+1}(\rho^{i}_q+\theta_{q+1}^{(c,i)}+\theta_{q+1}^{(o,i)}+\theta_{q+1}^{(g,i)}) \|_{L^1_tL^1}\notag\\
&\leq \|v_q\|_{C_{t,x}^0}\|\theta_{q+1}^{i}\|_{C_tL^1}+(\|\rho^{i}_q\|_{C_{t,x}^0}+\|\theta_{q+1}^{(c,i)}\|_{C_{t,x}^0}
+\|\theta_{q+1}^{(o,i)}\|_{C_{t,x}^0}+\|\theta_{q+1}^{(g,i)}\|_{C_{t,x}^0}
)\|w_{q+1}\|_{L^r_tL^p}\notag\\
&\lesssim (C_0\lambda_{q}^{d+4}+1)\lambda_{q+1}^{(12d+40)\alpha-\frac1{2N}}
\lesssim  C_0 \lambda_{q+1}^{(12d+41)\alpha-\frac{1}{2N}}\lesssim   C_0 \lambda_{q+1}^{-\alpha},\notag
\end{align}
where we used  \eqref{para42} and   conditions on the parameters to have $(12d+42)\alpha<\frac{1}{2N}$.

Since $w_{q+1}(t)=0$ on $[0,T_{q+1}]$, by  \eqref{5:bd:rhoin} and \eqref{bd:4wq+1lpr}    we obtain 
\begin{align}
\|w_{q+1} \rho^{in,i}\|_{L^1_tL^1}\lesssim \| \rho^{in,i}\|_{C_{[\frac1{6},1],x}^0}\|w_{q+1}\|_{L^r_tL^p}\lesssim C_{in} \lambda_{q+1}^{-\alpha}.\notag
\end{align} 

By  the estimates in \eqref{bd:4wq+1clp} and \eqref{bd:4theq+1pl2}   we have
\begin{align}
   \|w_{q+1}^{(c,i)} \theta_{q+1}^{(p,i)}\|_{L^1_tL^1}&\leq    \| \theta_{q+1}^{(p,i)}\|_{L^{d_0'}_tL^{d_0'}}\|w_{q+1}^{(c,i)}\|_{L^{d_0}_tL^{d_0}}\lesssim C_\rho C_0^{1/d_0'} l^{-3d-8}\frac{r_\perp}{r_\parallel}\lesssim  C_0 \lambda_{q+1}^{(6d+16)\alpha-\frac1N} \lesssim   C_0  \lambda_{q+1}^{-\alpha},\notag
\end{align}
where we used  \eqref{para42} and  conditions on the parameters to have $(6d+17)\alpha<\frac{1}{N}$. We choose $a$ large enough to absorb the universal constant.

In summary, we have
\begin{align}\label{bd:mlinl1}
    \|M^{i}_{lin}\|_{L_t^1L^1}\lesssim C_0\lambda_{q+1}^{-\alpha}\leq \frac13C_0 \delta_{q+2},
\end{align}
where we choose $a$ large enough to absorb the constant.

Moreover, by \eqref{bd:4rhoqc1}, \eqref{bs:theq+1locc2}   and all the H\"older's estimate in Section \ref{sec:4estwq+1}-\ref{sec:4esttheq+1g}, we have
\begin{align}
   \|\partial_t &M^{i}_{lin}\|_{L_t^1L^1}+  \| \nabla M^{i}_{lin}\|_{L_t^1L^1}\notag\\ 
   &\lesssim \|\nabla\theta^{(loc,i)}_{q+1}\|_{C_{t,x}^1}+\|v_q\theta^{i}_{q+1}+w_{q+1}(\rho^{i}_q+\theta_{q+1}^{(c,i)}+\theta_{q+1}^{(o,i)}+\theta_{q+1}^{(g,i)}+\rho^{in,i}) +w_{q+1}^{(c,i)}\theta_{q+1}^{(p,i)}\|_{C_{t,x}^1}\notag\\ 
   &\lesssim \lambda_{q+1}^{d+5}+\lambda_{q+1}^{(24d+80)\alpha+2d+7}\leq \frac13\lambda_{q+1}^{2d+8},\label{bd:mlinc1}
\end{align}
where we used $(12d+40)\alpha<\frac13$, and chose $a$ large to abosrb the universal constant.

 \subsubsection{Estimate of commutator error $M^{i}_{com}$.}
By  \eqref{bd:4rql1} and \eqref{para42}  we obtain
 \begin{align*}
  \|M_{ com}^i\|_{L_t^1L^1}\lesssim l(\|\partial_t M^{i}_{q}\|_{L_t^1L^1}+  \| \nabla M^{i}_{q}\|_{L_t^1L^1})\lesssim C_0 l\lambda_q^{2d+8}\lesssim   C_0  \lambda_{q+1}^{-\alpha}\leq \frac13C_0 \delta_{q+2},\notag\\  \|\partial_t M^{i}_{com}\|_{L_t^1L^1}+  \| \nabla M^{i}_{com}\|_{L_t^1L^1}\lesssim \|\partial_t M^{i}_q\|_{L_t^1L^1}+  \| \nabla M^{i}_{q}\|_{L_t^1L^1}\lesssim \lambda_{q}^{2d+8}\leq \frac13\lambda_{q+1}^{2d+8},
   \end{align*} 
   which together with \eqref{bd:moscl1}-\eqref{bd:mlinc1} implies that
\begin{align*}
    \|M^{i}_{q+1} \|_{L^1_tL^1}\leq C_0 \delta_{q+2},\ \
\|\partial_t M^{i}_{q+1}\|_{L_t^1L^1}+  \| \nabla M^{i}_{q+1}\|_{L_t^1L^1}&\leq \lambda_{q+1}^{2d+8}.
\end{align*}

We finish the proof of Proposition \ref{prop:case4}.
 \\
            
\noindent{\bf Acknowledgment.} The authors gratefully acknowledge Prof. Xiangchan Zhu and Dr. Zimo Hao for their valuable suggestions and insightful discussions.

\appendix
 \renewcommand{\appendixname}{Appendix~\Alph{section}}
  \renewcommand{\theequation}{A.\arabic{equation}}
  \section{An overview of convex integration on the torus}\label{app:a}
In this section, we present an overview of how the convex integration method applies to the continuity equation on the torus. 
For detailed derivations and estimates, we refer the reader to~\cite{BCDL21}, 
and to~\cite{LRZ25} for the advection-diffusion case on the torus. 
Our goal here is to outline the main ideas and the structure of the construction rather than to reproduce the full technical details. 
We also add remarks on the additional difficulties that arise when extending the argument to the whole space $\mathbb{R}^d$. 
The new ideas required to overcome these issues are presented in detail in Section~\ref{sec:proofmaireuslt}.

We consider the following continuity equation on $\mathbb{T}^d$, with $d \geq 2$:
\begin{align}
    \partial_t\rho+\div (v\rho)&=0,\ \div v=0,\notag\\
    \rho(0)&=1.\label{eq:tpe}
\end{align}
By the divergence-free condition on the drift $v$, equation~\eqref{eq:tpe} admits the constant solution 
$\overline{\rho}(t) = 1$ which is  a probability density, w.r.t. the normed volume measures on $\mathbb{T}^d$. 
Our goal is therefore to construct a divergence-free drift~$v$ such that~\eqref{eq:tpe} 
admits another positive, non-constant solution.

The construction proceeds via induction indexed by $q\in\mathbb{N}$. At each step $q\in\mN_0$, we construct a pair $(v_q, \rho_q,  M_q)$  satisfying the following system:
\begin{align}
\partial_t \rho_q+\div(v_q \rho_q)&=-\div M_q, \ \div v_q=0,\label{eq:idea:rhov}
\end{align}
Here $M_q$ denotes a vector field. As $q \to \infty$, we  aim to prove that in some topology,
$M_q \to 0$ and $(v_q, \rho_q) \to (v,\rho),$
which is a weak solution of the transport equation. 
At the same time, we require that $\int_{\mathbb{T}^d}\rho_q=1$ and $\rho_q\geq0,$ so that the solution can be interpreted as a probability density.

At each iterative step, we construct perturbations 
$w_{q+1}=v_{q+1}-v_q,\theta_{q+1}=\rho_{q+1}-\rho_q$
such that $(v_{q+1}, \rho_{q+1})$ is a solution of \eqref{eq:tpe} at level $q+1$ 
with a smaller residual stress term $M_{q+1}$, which can be written as
 \begin{align*}
     -\div M_{q+1}=\partial_t\theta_{q+1}+\div(v_q\theta_{q+1}+w_{q+1}\rho_q)+\div (w_{q+1}\theta_{q+1}-M_q).
 \end{align*}
 The stress term $M_{q+1}$ is then defined by solving the above divergence equation 
with the aid of a linear differential operator of order $-1$, denoted by $\div^{-1}$. 
We remark that the inverse divergence of a periodic function remains periodic. 
However, on the whole space, handling the inverse operator $\div^{-1}$ becomes substantially more delicate.

 The construction of the perturbations forms the core of the convex integration scheme.
 We define the perturbations $(w_{q+1}, \theta_{q+1})$ as a sum of highly oscillatory building blocks 
in order to achieve a cancellation between the low-frequency part of the quadratic term $w_{q+1} \theta_{q+1}$ and the stress $M_q$. 
Roughly speaking, we introduce the principal part of the perturbation 
$(w_{q+1}^{(p)}, \theta_{q+1}^{(p)})$, which takes the form
 
  \begin{align*}
      w_{q+1}^{(p)}=\sum_{\xi}a_{\xi}(M_q)W_{\xi}(\lambda_{q+1}x),\ \  \theta^{(p)}_{q+1}=\sum_{\xi}a_{\xi}(M_q)\Theta_{\xi}(\lambda_{q+1}x),
  \end{align*}
  where $\Theta_{\xi}$ is a carefully chosen positive solution to the continuity equation transported by $W_{\xi}$. 
These building blocks are also oscillating at a high frequency $\lambda_{q+1}\in\mathbb{N}$. 
The amplitude coefficients $a_{\xi}$ are chosen such that
  \begin{align*}
      \sum_\xi a^2_{\xi}(M_q)\cdot\int_{\mathbb{T}^d}( W_{\xi}\Theta_{\xi})(\lambda_{q+1}x)\dif x = \sum_\xi a^2_{\xi}(M_q)\xi=M_q.
  \end{align*}
  Here, the existence of such coefficients $a_{\xi}$ is ensured by the Geometric Lemma \ref{lem:cv2}. 
Roughly speaking, the low-frequency part of the product $w_{q+1}^{(p)} \theta_{q+1}^{(p)}$ cancels the stress term $M_q$. 
Additionally, to ensure that $w_{q+1}$ is divergence-free, we introduce a divergence-free corrector $w_{q+1}^{(c)}$ 
such that 
$$
w_{q+1} := w_{q+1}^{(p)} + w_{q+1}^{(c)}
$$
is divergence-free. Similarly, to guarantee that 
$
\int_{\mathbb{T}^d} \rho_{q+1} \, dx - \int_{\mathbb{T}^d} \rho_q \, dx = 0,
$
the actual perturbation for $\rho_{q+1}$ is defined as
  $$\theta_{q+1}:=\theta_{q+1}^{(p)}-\int_{\mathbb{T}^d}\theta_{q+1}^{(p)}\dif x.$$
  
Moreover, by taking advantage of the fact that the constructed $\theta^{(p)}_{q+1}$ is a positive function,   we have 
  \begin{align*}
      \rho_{q+1}-\rho_q\geq -\int_{\mathbb{T}^d}\theta^{(p)}_{q+1}\dif x ,
  \end{align*}
  which implies that the limited function $\rho$ satisfies $\rho\geq\rho_0+\sum_q ( \rho_{q+1}-\rho_q)\geq \inf_x\rho_0 -\sum_q\|\theta_{q+1}\|_{L^1}>0$ by choosing suitable parameters. Here, we use the compactness of the domain to ensure that $\rho_0$ is bounded below, 
so that the solution remains positive by keeping the perturbations $\theta_{q+1}$ sufficiently small. 
In contrast, on the whole space $\mathbb{R}^d$, this argument no longer applies, 
since the density $\rho_0$ decays to zero at infinity. 
Therefore, in the whole-space setting, the perturbations must be estimated carefully 
to ensure compatibility with the decay of $\rho_0$.  

 It then suffices to verify that $M_{q+1}$ is small in the Lebesgue space $L^1$. 
We choose the frequencies super-geometrically, with $\lambda_{q+1} \gg \lambda_q$ for $q \in \mathbb{N}$. 
The function 
$
v_q = \sum_{i \leq q} (v_i - v_{i-1}) + v_0 = \sum_{i \leq q} w_i + v_0
$
is thus a sum of perturbations of frequencies higher than $\lambda_q$, 
which are considered as low-frequencies compared with the building blocks $W_\xi(\lambda_{q+1} x)$ and $\Theta_\xi(\lambda_{q+1} x)$. 
Similarly, the functions $\rho_q$ and $M_q$ are also sums of perturbations with frequencies exceeding $\lambda_q$.
Then, we have 

  \begin{align*}
      \div (w^{(p)}_{q+1}\theta^{(p)}_{q+1}-M_q)&\sim \div\(\sum_\xi a_{\xi}^2(M_q)\mathbb{P}_{\geq\lambda_{q+1}/2}[(W_{\xi}\Theta_{\xi})(\lambda_{q+1}x)]\)\\
      &\sim \sum_\xi (\nabla a_{\xi}^2(M_q))\cdot\mathbb{P}_{\geq\lambda_{q+1}/2}[(W_{\xi}\Theta_{\xi})(\lambda_{q+1}x)],
  \end{align*}
 and as a consequence,
  \begin{align*}
\|\div^{-1}\div(w_{q+1}\theta_{q+1}-M_q)\|_{L^1}\lesssim \frac{1}{\lambda_{q+1}}\sum_\xi\|a^2_{\xi}(M_q)\|_{C^1}\|\mathbb{P}_{\geq\lambda_{q+1}/2}(W_{\xi}\Theta_{\xi})\|_{L^1}\lesssim \frac{\sum_\xi\|a^2_{\xi}(M_q)\|_{C^1}}{\lambda_{q+1}}.
\end{align*}
Here, the factor $\lambda_{q+1}^{-1}$ arises from the $-1$ order inverse divergence operator $\div^{-1}$. 
In the last term, the amplitude $a_\xi$ is of low frequency, of order at most $\lambda_q$, 
and is therefore very small due to the super-geometric growth $\lambda_{q+1} \gg \lambda_q$. 

Here, we note that on the whole space, the stress term must be assumed to have compact support 
in order to fit within the convex integration framework. 
However, in this setting, one cannot directly apply the inverse divergence operator, 
since even for a function with compact support, the inverse of its divergence 
may still have support extending over the entire space. In this paper, we introduce a new decomposition of the stress term   that separates the stress term into two distinct components.

\end{document}